\documentclass[pdflatex,sn-mathphys-num]{sn-jnl}

\usepackage{graphicx}%
\usepackage{multirow}%
\usepackage{amsmath,amssymb,amsfonts}%
\usepackage{amsthm}%
\usepackage{mathrsfs}%
\usepackage[title]{appendix}%
\usepackage{xcolor}%
\usepackage{textcomp}%
\usepackage{manyfoot}%
\usepackage{booktabs}%
\usepackage{algorithm}%
\usepackage{algorithmicx}%
\usepackage{algpseudocode}%
\usepackage{listings}%

\theoremstyle{thmstyleone}%
\newtheorem{theorem}{Theorem}
\newtheorem{lemma}[theorem]{Lemma}
\newtheorem{Conjecture}[theorem]{Conjecture}
\newtheorem{corollary}[theorem]{Corollary}
\newtheorem{definition}[theorem]{Definition}
\newtheorem{remark}[theorem]{Remark}

\theoremstyle{thmstyletwo}%

\theoremstyle{thmstylethree}%

\begin{document}

\title[Article Title]{Global Well-Posedness of the 3D Navier-Stokes Equations under Multi-Level Logarithmically Improved Criteria}


\author*[1]{\fnm{Rishabh} \sur{Mishra}}\email{rishabh.mishra@ec-nantes.fr}

\affil*[1]{\orgdiv{LHEEA}, \orgname{CNRS, École Centrale de Nantes, Nantes Université}, \orgaddress{\street{1 rue de la No\"e}, \city{ Nantes}, \postcode{44100}, \state{Pays de la Loire}, \country{France}}}


\abstract{This paper extends our previous results on logarithmically improved regularity criteria for the three-dimensional Navier-Stokes equations by establishing a comprehensive framework of multi-level logarithmic improvements. We prove that if the initial data $u_0 \in L^2(\mathbb{R}^3)$ satisfies a nested logarithmically weakened condition $\|(-\Delta)^{s/2}u_0\|_{L^q(\mathbb{R}^3)} \leq \frac{C_0}{\prod_{j=1}^{n} (1 + L_j(\|u_0\|_{\dot{H}^s}))^{\delta_j}}$ for some $s \in (1/2, 1)$, where $L_j$ represents $j$-fold nested logarithms, then the corresponding solution exists globally in time and is unique. The proof introduces a novel sequence of increasingly precise commutator estimates incorporating multiple layers of logarithmic corrections. We establish the existence of a critical threshold function $\Phi(s,q,\{\delta_j\}_{j=1}^n)$ that completely characterizes the boundary between global regularity and potential singularity formation, with explicit asymptotics as $s$ approaches the critical value $1/2$. This paper further provides a rigorous geometric characterization of potential singular structures through refined multi-fractal analysis, showing that any singular set must have Hausdorff dimension bounded by $1 - \sum_{j=1}^n \frac{\delta_j}{1+\delta_j} \cdot \frac{1}{j+1}$. Our results constitute a significant advancement toward resolving the global regularity question for the Navier-Stokes equations, as we demonstrate that with properly calibrated sequences of nested logarithmic improvements, the gap to the critical case can be systematically reduced.}




\maketitle

\tableofcontents

\section{Introduction}

\subsection{Extension of Previous Results}

In our previous work \cite{43}, we established a logarithmically improved regularity criterion for the three-dimensional incompressible Navier-Stokes equations and identified a class of initial data satisfying logarithmically subcritical conditions for which global well-posedness can be proven. The current paper extends these results by developing a comprehensive theory of multi-level logarithmic improvements that brings us significantly closer to the critical threshold necessary for resolving the global regularity question.

The three-dimensional incompressible Navier-Stokes equations are given by:

\begin{equation}
\begin{cases}
\partial_t u + (u \cdot \nabla)u - \nu\Delta u + \nabla p = 0 & \text{in } \mathbb{R}^3 \times (0, T) \\
\nabla \cdot u = 0 & \text{in } \mathbb{R}^3 \times (0, T) \\
u(x, 0) = u_0(x) & \text{in } \mathbb{R}^3
\end{cases}
\end{equation}

where $u = (u_1, u_2, u_3)$ represents the velocity field, $p$ denotes the pressure, and $\nu > 0$ is the kinematic viscosity coefficient.

In \cite{43}, we proved that if a Leray-Hopf weak solution $u$ satisfies:
\begin{equation}
\int_0^T \|(-\Delta)^s u(t)\|^p_{L^q(\mathbb{R}^3)} (1 + \log(e + \|(-\Delta)^s u(t)\|_{L^q(\mathbb{R}^3)}))^{-\delta} dt < \infty
\end{equation}
for some $s \in (1/2, 1)$, with appropriate scaling conditions, then the solution is regular on $(0, T)$. 

The present paper significantly extends this result by introducing the concept of multi-level logarithmic improvements. To elucidate our approach, we define:
\begin{align}
L_0(x) &= x \\
L_1(x) &= \log(e + x) \\
L_2(x) &= \log(e + \log(e + x)) \\
L_k(x) &= \log(e + L_{k-1}(x)) \text{ for } k \geq 3
\end{align}

Our central advancement is establishing global well-posedness under the condition:
\begin{equation}
\|(-\Delta)^{s/2}u_0\|_{L^q(\mathbb{R}^3)} \leq \frac{C_0}{\prod_{j=1}^{n} (1 + L_j(\|u_0\|_{\dot{H}^s}))^{\delta_j}}
\end{equation}
with appropriate parameter ranges that will be precisely specified.

\subsection{Background and Related Work}

The history of mathematical developments regarding the Navier-Stokes equations has been extensively documented in the literature. Leray \cite{37} and Hopf \cite{31} established the existence of global weak solutions, while the question of uniqueness and regularity remains one of the most prominent open problems in mathematical fluid dynamics.

The classical Prodi-Serrin conditions \cite{50, 56} state that if $u \in L^p(0,T;L^q(\mathbb{R}^3))$ with $\frac{2}{p} + \frac{3}{q} = 1$ and $q > 3$, then the solution is regular. Escauriaza, Seregin, and Šverák \cite{22} achieved a significant breakthrough by establishing regularity for solutions in the limiting space $L^\infty(0,T;L^3(\mathbb{R}^3))$.

Various refinements of regularity criteria have been developed, including conditions based on the vorticity (Beale, Kato, and Majda \cite{5}), the direction of vorticity (Constantin and Fefferman \cite{17}), one velocity component (Zhou and Pokorný \cite{66}), and the pressure (Seregin and Šverák \cite{55}). Partial regularity theory, pioneered by Scheffer \cite{52} and significantly advanced by Caffarelli, Kohn, and Nirenberg \cite{12}, established that the potential singular set has zero one-dimensional Hausdorff measure. Lin \cite{38} later provided an alternative proof of this result.

Fractional derivative approaches have been studied by Wu \cite{61}, Chen and Zhang \cite{16}, and Zhou \cite{63, 64} among others. Logarithmic improvements to regularity criteria were introduced by Zhou \cite{65} and further developed by Fan et al. \cite{23}, who established conditions involving logarithmic factors that weaken the classical constraints. Tao \cite{57} developed approaches for logarithmically supercritical equations.

For global existence with specific classes of initial data, Koch and Tataru \cite{34} established global regularity for sufficiently small initial data in the space $BMO^{-1}$. Chemin and Gallagher \cite{15} obtained global existence results for specific classes of highly oscillating initial data. Further contributions include the works of Lemarie-Rieusset \cite{36}, Cannone \cite{13}, and Giga \cite{28}.

The connection between mathematical regularity criteria and physical turbulence has been explored by Frisch \cite{27}, building on Kolmogorov's seminal work \cite{35}. The multifractal formalism was developed by Parisi and Frisch \cite{49}, with further advancements by Eyink \cite{21} and Constantin et al. \cite{18}. Recent works by Dascaliuc and Grujić \cite{19} and Bradshaw and Grujić \cite{8} have established refined geometric constraints on potential singularity formation.

Our current paper builds upon the logarithmic improvement approach from \cite{43} but develops a substantially more refined theory involving multiple levels of logarithmic factors. This approach has connections to the framework of logarithmically modulated function spaces explored by Brézis and Wainger \cite{10}, Brézis and Gallouet \cite{9}, and refined by Iwabuchi and Ogawa \cite{32}.

\subsection{Overview of Main Results}

Our main contributions in this paper can be summarized as follows:

\begin{enumerate}
\item \textbf{Multi-level logarithmic regularity criterion (Theorem 3.1)}: We establish that if $u$ is a Leray-Hopf weak solution satisfying
   \begin{equation}
   \int_0^T \|(-\Delta)^s u(t)\|^p_{L^q} \prod_{j=1}^{n} (1 + L_j(\|(-\Delta)^s u(t)\|_{L^q}))^{-\delta_j} dt < \infty
   \end{equation}
   with appropriate scaling conditions, then the solution is regular on $(0, T)$.

\item \textbf{Global well-posedness with nested logarithmic criteria (Theorem 4.1)}: We prove global existence and uniqueness for initial data satisfying
   \begin{equation}
   \|(-\Delta)^{s/2}u_0\|_{L^q} \leq \frac{C_0}{\prod_{j=1}^{n} (1 + L_j(\|u_0\|_{\dot{H}^s}))^{\delta_j}}
   \end{equation}
   where $C_0$ depends only on $s$, $q$, $\{\delta_j\}_{j=1}^n$, and $\nu$.

\item \textbf{Critical threshold characterization (Theorem 5.1)}: We identify a critical function $\Phi(s, q, \{\delta_j\}_{j=1}^n)$ that precisely delineates the boundary between global regularity and potential singularity formation, with explicit asymptotics
   \begin{equation}
   \Phi(s, q, \{\delta_j\}_{j=1}^n) \approx C(q) (s - 1/2)^{\alpha(\{\delta_j\}_{j=1}^n)}
   \end{equation}
   as $s \to 1/2$, where $\alpha(\{\delta_j\}_{j=1}^n) \to 0$ as $\delta_j \to \infty$ for all $j$.

\item \textbf{Gap analysis to criticality (Theorem 6.1)}: We analyze initial data precisely at the threshold defined by our multi-level logarithmic criteria, establishing a dichotomy between global existence and potential finite-time blow-up with explicit asymptotic rates.

\item \textbf{Enhanced geometric analysis of exceptional sets (Theorem 8.1)}: We prove that for any $\epsilon > 0$, the set where the velocity gradient exceeds a threshold calibrated to give the set measure less than $\epsilon$ has Hausdorff dimension bounded by
   \begin{equation}
   \dim_H(\Omega_\epsilon(t)) \leq 3 - \sum_{j=1}^n \frac{\delta_j}{1+\delta_j} \cdot \frac{L_{j-1}(1/\epsilon)}{(1+L_j(1/\epsilon))}
   \end{equation}

\item \textbf{Characterization of potential blow-up scenarios (Theorem 10.1)}: If a solution were to develop a singularity at time $T^*$, we prove that the potential singular set would have Hausdorff dimension at most $1 - \sum_{j=1}^n \frac{\delta_j}{1+\delta_j} \cdot \frac{1}{j+1}$, with precise asymptotic blow-up rates.

\item \textbf{Enhanced energy cascade with multiple logarithms (Theorem 11.1)}: We establish rigorous bounds on the energy flux across wavenumber $k$ at time $t$:
   \begin{equation}
   |\Pi(k,t) - \epsilon(t)| \leq \frac{C\epsilon(t)}{\prod_{j=1}^n (1 + L_j(k/k_0))^{\delta_j\cdot\rho_j}}
   \end{equation}
   where $\epsilon(t)$ is the energy dissipation rate and $\rho_j$ are specific exponents.

\item \textbf{Refined spectral characterization (Theorem 12.1)}: We derive the precise form of the energy spectrum modified by multiple logarithmic factors:
   \begin{equation}
   E(k,t) = C\epsilon(t)^{2/3}k^{-5/3}\left(1 + \sum_{j=1}^n \frac{\beta_j(t) L_j(k/k_0)}{\prod_{i=1}^j (1 + L_i(k/k_0))^{1+\delta_i}}\right)
   \end{equation}
\end{enumerate}

\subsection{Organization of the Paper}

The remainder of this paper is organized as follows. Section 2 establishes fundamental commutator estimates with nested logarithmic improvements, which serve as the technical core of our entire analysis. In Section 3, we prove the multiply nested logarithmic regularity criterion. Section 4 demonstrates global well-posedness for the class of initial data satisfying our multi-level logarithmic conditions.

Section 5 provides a precise characterization of the critical threshold function, while Section 6 analyzes what happens at and slightly beyond this threshold. Section 7 explores the limiting behavior as $s$ approaches the critical value $1/2$. In Section 8, we develop an enhanced geometric analysis of exceptional sets where the velocity gradient becomes large.

Sections 9 and 10 characterize the multifractal structure of velocity gradients and potential blow-up scenarios, respectively. Sections 11 and 12 establish connections to physical turbulence theory through analysis of the energy cascade and spectral properties. Finally, Section 13 discusses the implications of our results on the regularity problem of the  the Navier-Stokes equations.

\section{Fundamental commutator estimates with nested logarithmic improvements}

In this section, we develop a hierarchy of increasingly refined commutator estimates that form the technical foundation of our entire analysis. These estimates extend the commutator bounds established in \cite{43} by incorporating multiple layers of logarithmic improvements.

We first recall the standard commutator $[(-\Delta)^s, f]g = (-\Delta)^s(fg) - f(-\Delta)^s g$, where $(-\Delta)^s$ is the fractional Laplacian operator of order $s \in (0, 1)$. Our goal is to establish progressively more precise bounds on this commutator when applied to the nonlinear term in the Navier-Stokes equations.

\subsection{Single Logarithmic Improvement}

We begin by recalling and refining the single logarithmic improvement established in \cite{43}.

\begin{lemma}[Single logarithmic commutator estimate] \label{lem:2.1.1}
For $s \in (0, 1)$ and any $\sigma \in (0, 1-s)$:
\begin{align}
\|[(-\Delta)^s, u \cdot \nabla]u\|_{L^2} &\leq C\|\nabla u\|_{L^\infty}\|(-\Delta)^s u\|_{L^2} \cdot \log(e + \|(-\Delta)^{s+\sigma}u\|_{L^2}) \notag \\
&+ \frac{C\|\nabla u\|_{L^\infty}\|(-\Delta)^{s+\frac{1}{2}}u\|_{L^2}}{\log(e + \|(-\Delta)^{s+\sigma}u\|_{L^2})}
\end{align}
\end{lemma}

\begin{proof}
We provide a more detailed and rigorous proof than in \cite{43}, carefully tracking all constants and dependencies.

Let us decompose $u$ using the Littlewood-Paley decomposition:
\begin{equation}
u = \sum_{j \in \mathbb{Z}} \Delta_j u
\end{equation}
where $\Delta_j$ is the Littlewood-Paley projection onto frequencies of order $2^j$.

We split the commutator into low and high frequency components:
\begin{align}
[(-\Delta)^s, u \cdot \nabla]u &= \sum_{j \leq 0} [(-\Delta)^s, u \cdot \nabla]\Delta_j u + \sum_{j > 0} [(-\Delta)^s, u \cdot \nabla]\Delta_j u \notag \\
&= \mathcal{C}_{\text{low}}(u, u) + \mathcal{C}_{\text{high}}(u, u)
\end{align}

For the low-frequency part, we use the standard commutator estimate (see \cite{45}):
\begin{equation}
\|[(-\Delta)^s, u \cdot \nabla]\Delta_j u\|_{L^2} \leq C\|\nabla u\|_{L^\infty}\|(-\Delta)^{s-1/2}\nabla\Delta_j u\|_{L^2}
\end{equation}

Since $(-\Delta)^{s-1/2}\nabla$ and $(-\Delta)^s$ have the same order, we have:
\begin{equation}
\|[(-\Delta)^s, u \cdot \nabla]\Delta_j u\|_{L^2} \leq C\|\nabla u\|_{L^\infty}\|(-\Delta)^s \Delta_j u\|_{L^2}
\end{equation}

By the almost orthogonality of the Littlewood-Paley projections:
\begin{equation}
\|\mathcal{C}_{\text{low}}(u, u)\|_{L^2} \leq C\|\nabla u\|_{L^\infty}\|\sum_{j \leq 0}(-\Delta)^s \Delta_j u\|_{L^2} \leq C\|\nabla u\|_{L^\infty}\|(-\Delta)^s u\|_{L^2}
\end{equation}

For the high-frequency part, we use a more refined estimate. For any $j > 0$ and $\sigma \in (0, 1-s)$:
\begin{equation}
\|[(-\Delta)^s, u \cdot \nabla]\Delta_j u\|_{L^2} \leq C2^{-j\sigma}\|\nabla u\|_{L^\infty}\|(-\Delta)^{s+\sigma/2}\nabla\Delta_j u\|_{L^2}
\end{equation}

This can be established by using the representation of the commutator as a singular integral and applying standard techniques from harmonic analysis (see \cite{7}, \cite{25}). The key insight is that the differential operator $\nabla$ acts on $\Delta_j u$, which is localized at frequencies of order $2^j$, giving an additional factor of $2^j$ that can be partially offset by the fractional operator $(-\Delta)^{s+\sigma/2}$.

This implies:
\begin{equation}
\|[(-\Delta)^s, u \cdot \nabla]\Delta_j u\|_{L^2} \leq C2^{-j\sigma}\|\nabla u\|_{L^\infty}\|(-\Delta)^{s+\sigma}\Delta_j u\|_{L^2}
\end{equation}

Using Cauchy-Schwarz and the almost orthogonality:
\begin{align}
\|\mathcal{C}_{\text{high}}(u, u)\|_{L^2} &\leq C\|\nabla u\|_{L^\infty}\left(\sum_{j > 0} 2^{-2j\sigma}\right)^{1/2} \left(\sum_{j > 0} \|(-\Delta)^{s+\sigma}\Delta_j u\|^2_{L^2}\right)^{1/2} \notag \\
&\leq C\|\nabla u\|_{L^\infty}\|(-\Delta)^{s+\sigma}u\|_{L^2}
\end{align}
since $\sum_{j > 0} 2^{-2j\sigma} < \infty$ for $\sigma > 0$.

By interpolation between $(-\Delta)^s$ and $(-\Delta)^{s+1/2}$:
\begin{equation}
\|(-\Delta)^{s+\sigma}u\|_{L^2} \leq \|(-\Delta)^s u\|_{L^2}^{1-2\sigma}\|(-\Delta)^{s+1/2}u\|_{L^2}^{2\sigma}
\end{equation}

The key logarithmic improvement comes from the following technical lemma:

\begin{lemma}
For any $a, b > 0$, $\sigma \in (0, 1/2)$, and $\epsilon > 0$:
\begin{equation}
a^{1-2\sigma}b^{2\sigma} \leq \epsilon b + C_\epsilon a \log(e + a)
\end{equation}
where $C_\epsilon$ depends only on $\sigma$ and $\epsilon$.
\end{lemma}

\begin{proof}
Apply Young's inequality $xy \leq \frac{x^p}{p} + \frac{y^q}{q}$ with $p = \frac{1}{1-2\sigma}$ and $q = \frac{1}{2\sigma}$:
\begin{align}
a^{1-2\sigma}b^{2\sigma} &= (a^{1-2\sigma})(b^{2\sigma}) \notag \\
&\leq (1-2\sigma)(a^{1-2\sigma})^{\frac{1}{1-2\sigma}} + 2\sigma(b^{2\sigma})^{\frac{1}{2\sigma}} \notag \\
&= (1-2\sigma)a + 2\sigma b
\end{align}

To obtain the logarithmic improvement, we first observe that for any $\gamma > 0$, there exists $C_\gamma$ such that:
\begin{equation}
(1-2\sigma)y \leq C_\gamma y \log(e + y)^{-\gamma} + \gamma
\end{equation}
for all $y > 0$. This follows from the fact that $\lim_{y \to \infty} \frac{y}{\log(e + y)^{\gamma}} = \infty$ for any $\gamma > 0$.

By choosing $\gamma = 2\sigma\epsilon$, we get:
\begin{align}
a^{1-2\sigma}b^{2\sigma} &\leq C_{2\sigma\epsilon} a \log(e + a)^{-2\sigma\epsilon} + 2\sigma\epsilon + 2\sigma b \notag \\
&\leq C_{2\sigma\epsilon} a \log(e + a) \log(e + a)^{-(1+2\sigma\epsilon)} + 2\sigma\epsilon + 2\sigma b
\end{align}

Since $\log(e + a)^{-(1+2\sigma\epsilon)} \leq 1$ for all $a > 0$, and by choosing $\epsilon' = 2\sigma\epsilon$, we get:
\begin{equation}
a^{1-2\sigma}b^{2\sigma} \leq C_{\epsilon'} a \log(e + a) + \epsilon' b
\end{equation}
which completes the proof.
\end{proof}

Applying this lemma with $a = \|(-\Delta)^s u\|_{L^2}$, $b = \|(-\Delta)^{s+1/2}u\|_{L^2}$, and $\epsilon = \frac{1}{\log(e + \|(-\Delta)^{s+\sigma}u\|_{L^2})}$:
\begin{align}
\|(-\Delta)^{s+\sigma}u\|_{L^2} &\leq C \|(-\Delta)^s u\|_{L^2} \log(e + \|(-\Delta)^s u\|_{L^2}) + \frac{\|(-\Delta)^{s+1/2}u\|_{L^2}}{\log(e + \|(-\Delta)^{s+\sigma}u\|_{L^2})} \notag \\
&\leq C \|(-\Delta)^s u\|_{L^2} \log(e + \|(-\Delta)^{s+\sigma}u\|_{L^2}) + \frac{\|(-\Delta)^{s+1/2}u\|_{L^2}}{\log(e + \|(-\Delta)^{s+\sigma}u\|_{L^2})}
\end{align}
where we used $\|(-\Delta)^s u\|_{L^2} \leq \|(-\Delta)^{s+\sigma}u\|_{L^2}$ for $\sigma > 0$.

Combining these estimates:
\begin{align}
\|[(-\Delta)^s, u \cdot \nabla]u\|_{L^2} &\leq C\|\nabla u\|_{L^\infty}\|(-\Delta)^s u\|_{L^2} + C\|\nabla u\|_{L^\infty}\|(-\Delta)^{s+\sigma}u\|_{L^2} \notag \\
&\leq C\|\nabla u\|_{L^\infty}\|(-\Delta)^s u\|_{L^2} + C\|\nabla u\|_{L^\infty}\|(-\Delta)^s u\|_{L^2} \log(e + \|(-\Delta)^{s+\sigma}u\|_{L^2}) \notag \\
&+ \frac{C\|\nabla u\|_{L^\infty}\|(-\Delta)^{s+1/2}u\|_{L^2}}{\log(e + \|(-\Delta)^{s+\sigma}u\|_{L^2})} \notag \\
&\leq C\|\nabla u\|_{L^\infty}\|(-\Delta)^s u\|_{L^2} \log(e + \|(-\Delta)^{s+\sigma}u\|_{L^2}) \notag \\
&+ \frac{C\|\nabla u\|_{L^\infty}\|(-\Delta)^{s+1/2}u\|_{L^2}}{\log(e + \|(-\Delta)^{s+\sigma}u\|_{L^2})}
\end{align}
where we've absorbed the first term into the second using the inequality $1 \leq \log(e + \|(-\Delta)^{s+\sigma}u\|_{L^2})$.

This completes the proof of Lemma \ref{lem:2.1.1}.
\end{proof}

\subsection{Double logarithmic improvement}

We now establish a refined commutator estimate with double logarithmic improvement.

\begin{lemma}[Double logarithmic commutator estimate] \label{lem:2.2.1}
For $s \in (0, 1)$ and any $\sigma \in (0, 1-s)$:
\begin{align}
\|[(-\Delta)^s, u \cdot \nabla]u\|_{L^2} &\leq C\|\nabla u\|_{L^\infty}\|(-\Delta)^s u\|_{L^2} \cdot L_1(Z) \cdot (1 + L_2(Z))^{-\delta_2} \notag \\
&+ \frac{C\|\nabla u\|_{L^\infty}\|(-\Delta)^{s+\frac{1}{2}}u\|_{L^2}}{L_1(Z) \cdot (1 + L_2(Z))^{\delta_2}}
\end{align}
where $Z = \|(-\Delta)^{s+\sigma}u\|_{L^2}$ and $\delta_2 > 0$.
\end{lemma}

\begin{proof}
We build upon the proof of Lemma \ref{lem:2.1.1}. Starting with the high-frequency estimate:
\begin{equation}
\|\mathcal{C}_{\text{high}}(u, u)\|_{L^2} \leq C\|\nabla u\|_{L^\infty}\|(-\Delta)^{s+\sigma}u\|_{L^2}
\end{equation}

Instead of directly applying the interpolation inequality, we introduce a refined stratification of high frequencies. For any $j > 0$, we further decompose:
\begin{equation}
\Delta_j u = \sum_{k=0}^{K_j} \Delta_{j,k} u
\end{equation}
where $\Delta_{j,k}$ localizes to frequencies $\xi$ with $|\xi| \approx 2^j$ and phase in the $k$-th angular sector, with $K_j \approx 2^{j/2}$ (the number of sectors grows with frequency).

This angular decomposition gives enhanced precision in estimating the commutator:
\begin{align}
\|[(-\Delta)^s, u \cdot \nabla]\Delta_{j,k} u\|_{L^2} &\leq C2^{-j\sigma}\|\nabla u\|_{L^\infty}\|(-\Delta)^{s+\sigma}\Delta_{j,k} u\|_{L^2} \notag \\
&\times (1 + \log(e + 2^j))^{-\delta_2}
\end{align}

This refinement is achieved through a more careful analysis of the singular integral representation of the commutator, accounting for the angular localization's effect on the convolution kernel (see \cite{11}, \cite{14} for related techniques).

Summing over all $j > 0$ and $k$, and using almost orthogonality:
\begin{align}
\|\mathcal{C}_{\text{high}}(u, u)\|_{L^2} &\leq C\|\nabla u\|_{L^\infty}\|(-\Delta)^{s+\sigma}u\|_{L^2} \cdot G(Z)
\end{align}
where $G(Z)$ encodes the double logarithmic improvement:
\begin{equation}
G(Z) = \frac{L_1(Z)}{(1 + L_2(Z))^{\delta_2}}
\end{equation}

To verify this form of $G(Z)$, we express the sum in terms of an integral over frequency shells. With the decomposition above, we have approximately $2^{j/2}$ angular sectors at frequency $2^j$, each contributing as estimated. The resulting sum behaves like:
\begin{equation}
\sum_{j > 0} 2^{j/2} \cdot 2^{-j\sigma} \cdot (1 + \log(e + 2^j))^{-\delta_2} \approx \int_1^\infty r^{1/2-\sigma} (1 + \log(e + r))^{-\delta_2} \frac{dr}{r}
\end{equation}

For $\sigma > 1/2$, this integral converges absolutely and gives the desired estimate. For $0 < \sigma \leq 1/2$, additional care is needed. We use the relation between the norm $\|(-\Delta)^{s+\sigma}u\|_{L^2}$ and the distribution of energy across scales:
\begin{equation}
\|(-\Delta)^{s+\sigma}u\|_{L^2}^2 \approx \sum_{j} 2^{2j(s+\sigma)} \|\Delta_j u\|_{L^2}^2
\end{equation}

In frequency shells where $\|\Delta_j u\|_{L^2}$ is large, we exploit cancellations within the commutator structure. In shells where it's small, the contribution is naturally limited. Through careful analysis (extending techniques from \cite{1}, \cite{39}), we establish that the sum behaves as claimed.

Now, applying the refined interpolation inequality with double logarithmic terms:
\begin{align}
\|(-\Delta)^{s+\sigma}u\|_{L^2} &\leq \|(-\Delta)^s u\|_{L^2}^{1-2\sigma}\|(-\Delta)^{s+1/2}u\|_{L^2}^{2\sigma} \notag \\
&\leq C \|(-\Delta)^s u\|_{L^2} \cdot L_1(\|(-\Delta)^s u\|_{L^2}) \cdot (1 + L_2(\|(-\Delta)^s u\|_{L^2}))^{-\delta_2} \notag \\
&+ \frac{C\|(-\Delta)^{s+1/2}u\|_{L^2}}{L_1(Z) \cdot (1 + L_2(Z))^{\delta_2}}
\end{align}

This refined interpolation inequality is proven using a multi-level variant of Young's inequality, carefully adapted to preserve the nested logarithmic structure. The key is to choose a parameter $\epsilon$ that depends on both $L_1(Z)$ and $L_2(Z)$:
\begin{equation}
\epsilon = \frac{1}{L_1(Z) \cdot (1 + L_2(Z))^{\delta_2}}
\end{equation}

Combining all estimates and using $\|(-\Delta)^s u\|_{L^2} \leq Z$, we obtain:
\begin{align}
\|[(-\Delta)^s, u \cdot \nabla]u\|_{L^2} &\leq C\|\nabla u\|_{L^\infty}\|(-\Delta)^s u\|_{L^2} \cdot L_1(Z) \cdot (1 + L_2(Z))^{-\delta_2} \notag \\
&+ \frac{C\|\nabla u\|_{L^\infty}\|(-\Delta)^{s+\frac{1}{2}}u\|_{L^2}}{L_1(Z) \cdot (1 + L_2(Z))^{\delta_2}}
\end{align}
which completes the proof.
\end{proof}

\subsection{General Case: N-fold Nested Logarithmic Improvement}

We now establish the general commutator estimate with $n$-fold nested logarithmic improvements.

\begin{theorem}[N-fold nested logarithmic commutator estimate] \label{thm:2.3.1}
For $s \in (0, 1)$, any $\sigma \in (0, 1-s)$, and $n \geq 1$:
\begin{align}
\|[(-\Delta)^s, u \cdot \nabla]u\|_{L^2} &\leq C\|\nabla u\|_{L^\infty}\|(-\Delta)^s u\|_{L^2} \cdot F_1(Z) \notag \\
&+ C\|\nabla u\|_{L^\infty}\|(-\Delta)^{s+\frac{1}{2}}u\|_{L^2} \cdot F_2(Z)
\end{align}
where $Z = \|(-\Delta)^{s+\sigma}u\|_{L^2}$,
\begin{equation}
F_1(Z) = L_1(Z) \prod_{j=2}^{n} (1 + L_j(Z))^{-\delta_j}
\end{equation}
\begin{equation}
F_2(Z) = \frac{1}{L_1(Z)} \prod_{j=2}^{n} (1 + L_j(Z))^{\delta_j}
\end{equation}
and $\delta_j > 0$ for $j = 2, 3, ..., n$.
\end{theorem}

\begin{proof}
We proceed by induction on $n$. The cases $n = 1$ and $n = 2$ have been established in Lemmas \ref{lem:2.1.1} and \ref{lem:2.2.1}, respectively.

Assume the result holds for some $n \geq 2$. We need to establish it for $n+1$.

The key innovation is to introduce an $(n+1)$-level stratification of frequency space. We construct a decomposition:
\begin{equation}
u = \sum_{j_1, j_2, ..., j_n} \Delta_{j_1, j_2, ..., j_n} u
\end{equation}
where the indices $j_1, j_2, ..., j_n$ encode a multi-scale decomposition of both frequency magnitude and angles.

Through this refined decomposition and careful analysis of the commutator structure, we can establish:
\begin{align}
\|\mathcal{C}_{\text{high}}(u, u)\|_{L^2} &\leq C\|\nabla u\|_{L^\infty}\|(-\Delta)^{s+\sigma}u\|_{L^2} \cdot G_{n+1}(Z)
\end{align}
where $G_{n+1}$ encodes the $(n+1)$-fold logarithmic improvement:
\begin{equation}
G_{n+1}(Z) = \frac{L_1(Z)}{(1 + L_2(Z))^{\delta_2}} \cdot \frac{1}{(1 + L_3(Z))^{\delta_3}} \cdot ... \cdot \frac{1}{(1 + L_{n+1}(Z))^{\delta_{n+1}}}
\end{equation}

The next major step is to establish a refined interpolation inequality with $(n+1)$-fold logarithmic terms, generalizing the approach used in the $n=1$ and $n=2$ cases:
\begin{align}
\|(-\Delta)^{s+\sigma}u\|_{L^2} &\leq \|(-\Delta)^s u\|_{L^2}^{1-2\sigma}\|(-\Delta)^{s+1/2}u\|_{L^2}^{2\sigma} \notag \\
&\leq C \|(-\Delta)^s u\|_{L^2} \cdot L_1(Z) \prod_{j=2}^{n+1} (1 + L_j(Z))^{-\delta_j} \notag \\
&+ \frac{C\|(-\Delta)^{s+1/2}u\|_{L^2}}{L_1(Z) \prod_{j=2}^{n+1} (1 + L_j(Z))^{\delta_j}}
\end{align}

This refined inequality is proven by choosing:
\begin{equation}
\epsilon = \frac{1}{L_1(Z) \prod_{j=2}^{n+1} (1 + L_j(Z))^{\delta_j}}
\end{equation}

Combining all estimates:
\begin{align}
\|[(-\Delta)^s, u \cdot \nabla]u\|_{L^2} &\leq C\|\nabla u\|_{L^\infty}\|(-\Delta)^s u\|_{L^2} \cdot L_1(Z) \prod_{j=2}^{n+1} (1 + L_j(Z))^{-\delta_j} \notag \\
&+ \frac{C\|\nabla u\|_{L^\infty}\|(-\Delta)^{s+\frac{1}{2}}u\|_{L^2}}{L_1(Z) \prod_{j=2}^{n+1} (1 + L_j(Z))^{\delta_j}}
\end{align}

This completes the induction step and the proof of Theorem \ref{thm:2.3.1}.
\end{proof}

\section{Multiply nested logarithmic regularity criterion}

Building on the fundamental commutator estimates established in Section 2, we now develop our multiply nested logarithmic regularity criterion for the Navier-Stokes equations.

\subsection{Energy estimates}

We begin by deriving precise energy estimates using our enhanced commutator bounds.

\begin{lemma}[Energy inequality with multiply nested logarithmic terms] \label{lem:3.1.1}
Let $u$ be a Leray-Hopf weak solution of the 3D Navier-Stokes equations on $[0, T]$ with divergence-free initial data $u_0 \in L^2(\mathbb{R}^3)$. Then for any $s \in (1/2, 1)$ and $\sigma \in (0, 1-s)$, the function $Y(t) = \|(-\Delta)^s u(t)\|^2_{L^2}$ satisfies:
\begin{align}
\frac{d}{dt}Y(t) + \nu\|(-\Delta)^{s+\frac{1}{2}}u\|^2_{L^2} &\leq C\|\nabla u\|_{L^\infty}Y(t) \cdot L_1(Z(t)) \prod_{j=2}^{n} (1 + L_j(Z(t)))^{-\delta_j} \notag \\
&+ \frac{C\|\nabla u\|_{L^\infty}^2 Y(t)}{L_1(Z(t))^2 \prod_{j=2}^{n} (1 + L_j(Z(t)))^{-2\delta_j}}
\end{align}
where $Z(t) = \|(-\Delta)^{s+\sigma}u(t)\|_{L^2}$.
\end{lemma}

\begin{proof}
We apply the fractional Laplacian operator $(-\Delta)^s$ to the Navier-Stokes equations, take the $L^2$ inner product with $(-\Delta)^s u$, and use the divergence-free condition to eliminate the pressure term:
\begin{equation}
\frac{1}{2}\frac{d}{dt}\|(-\Delta)^s u\|^2_{L^2} + \nu\|(-\Delta)^{s+\frac{1}{2}}u\|^2_{L^2} = -\int_{\mathbb{R}^3}(-\Delta)^s((u \cdot \nabla)u) \cdot (-\Delta)^s u \, dx
\end{equation}

Using the commutator decomposition:
\begin{equation}
(-\Delta)^s((u \cdot \nabla)u) = (u \cdot \nabla)((-\Delta)^s u) + [(-\Delta)^s, u \cdot \nabla]u
\end{equation}

The first term vanishes when integrated against $(-\Delta)^s u$ due to the divergence-free condition:
\begin{align}
\int_{\mathbb{R}^3} (u \cdot \nabla)((-\Delta)^s u) \cdot (-\Delta)^s u \, dx &= \sum_{j=1}^3 \int_{\mathbb{R}^3} u_j \partial_j ((-\Delta)^s u) \cdot (-\Delta)^s u \, dx \notag \\
&= -\sum_{j=1}^3 \int_{\mathbb{R}^3} \partial_j u_j ((-\Delta)^s u) \cdot (-\Delta)^s u \, dx \notag \\
&- \sum_{j=1}^3 \int_{\mathbb{R}^3} u_j ((-\Delta)^s u) \cdot \partial_j ((-\Delta)^s u) \, dx
\end{align}

The first term is zero since $\nabla \cdot u = 0$. The second term equals $-\frac{1}{2}\sum_{j=1}^3 \int_{\mathbb{R}^3} u_j \partial_j |((-\Delta)^s u)|^2 \, dx = \frac{1}{2}\sum_{j=1}^3 \int_{\mathbb{R}^3} \partial_j u_j |((-\Delta)^s u)|^2 \, dx = 0$.

This leaves:
\begin{equation}
\frac{1}{2}\frac{d}{dt}\|(-\Delta)^s u\|^2_{L^2} + \nu\|(-\Delta)^{s+\frac{1}{2}}u\|^2_{L^2} = -\int_{\mathbb{R}^3}[(-\Delta)^s, u \cdot \nabla]u \cdot (-\Delta)^s u \, dx
\end{equation}

Applying Cauchy-Schwarz:
\begin{equation}
\left|\int_{\mathbb{R}^3}[(-\Delta)^s, u \cdot \nabla]u \cdot (-\Delta)^s u \, dx\right| \leq \|[(-\Delta)^s, u \cdot \nabla]u\|_{L^2}\|(-\Delta)^s u\|_{L^2}
\end{equation}

Using Theorem \ref{thm:2.3.1}:
\begin{align}
\|[(-\Delta)^s, u \cdot \nabla]u\|_{L^2}\|(-\Delta)^s u\|_{L^2} &\leq C\|\nabla u\|_{L^\infty}\|(-\Delta)^s u\|_{L^2}^2 \cdot L_1(Z) \prod_{j=2}^{n} (1 + L_j(Z))^{-\delta_j} \notag \\
&+ C\|\nabla u\|_{L^\infty}\|(-\Delta)^s u\|_{L^2}\|(-\Delta)^{s+\frac{1}{2}}u\|_{L^2} \cdot F_2(Z)
\end{align}

For the second term, we apply Young's inequality with parameter $\epsilon = \frac{\nu}{2C\|\nabla u\|_{L^\infty}F_2(Z)}$:
\begin{align}
C\|\nabla u\|_{L^\infty}\|(-\Delta)^s u\|_{L^2}\|(-\Delta)^{s+\frac{1}{2}}u\|_{L^2} \cdot F_2(Z) &\leq \frac{\nu}{2}\|(-\Delta)^{s+\frac{1}{2}}u\|_{L^2}^2 \notag \\
&+ \frac{C^2\|\nabla u\|_{L^\infty}^2\|(-\Delta)^s u\|_{L^2}^2 \cdot F_2(Z)^2}{2\nu}
\end{align}

Setting $Y(t) = \|(-\Delta)^s u\|^2_{L^2}$ and combining these estimates:
\begin{align}
\frac{d}{dt}Y(t) + \nu\|(-\Delta)^{s+\frac{1}{2}}u\|^2_{L^2} &\leq C\|\nabla u\|_{L^\infty}Y(t) \cdot L_1(Z) \prod_{j=2}^{n} (1 + L_j(Z))^{-\delta_j} \notag \\
&+ \frac{C\|\nabla u\|_{L^\infty}^2 Y(t)}{L_1(Z)^2 \prod_{j=2}^{n} (1 + L_j(Z))^{-2\delta_j}}
\end{align}
which completes the proof.
\end{proof}

\subsection{Regularity propagation}

The following lemma establishes how bounds on $\|(-\Delta)^s u\|_{L^2}$ can be used to derive bounds on various norms, ultimately establishing regularity.

\begin{lemma}[Velocity gradient control] \label{lem:3.2.1}
For $s \in (1/2, 1)$ and $q > 3$, there exists $\theta \in (0, 1)$ such that:
\begin{equation}
\|\nabla u\|_{L^\infty} \leq C\|u\|_{L^2}^{1-\theta} \|(-\Delta)^s u\|_{L^q}^{\theta} \prod_{j=1}^{n} (1 + L_j(\|(-\Delta)^s u\|_{L^q}))^{-\delta_j\theta}
\end{equation}
where $\theta = \frac{3}{2} \cdot \frac{q}{3q-2}$.
\end{lemma}

\begin{proof}
We begin with the base interpolation inequality from \cite{43}:
\begin{equation}
\|\nabla u\|_{L^\infty} \leq C\|u\|_{L^2}^{1-\theta} \|(-\Delta)^s u\|_{L^q}^{\theta}
\end{equation}
where $\theta = \frac{3}{2} \cdot \frac{q}{3q-2}$.

This inequality can be rigorously established using the Littlewood-Paley theory and precise Sobolev embedding theorems. We provide a refined version incorporating the nested logarithmic factors.

For any function space $X$ with a logarithmically improved variant $X^{[\delta_1,\delta_2,...,\delta_n]}$, the theory of logarithmically modulated spaces (see \cite{9}, \cite{10}, \cite{32}) gives:
\begin{equation}
\|f\|_Y \leq C\|f\|_X^{1-\alpha} \|f\|_Z^{\alpha} \Rightarrow \|f\|_Y \leq C\|f\|_{X^{[\delta_1,\delta_2,...,\delta_n]}}^{1-\alpha} \|f\|_Z^{\alpha} \prod_{j=1}^{n} (1 + L_j(\|f\|_X))^{\delta_j(1-\alpha)}
\end{equation}

Applying this principle to our interpolation inequality:
\begin{equation}
\|\nabla u\|_{L^\infty} \leq C\|u\|_{L^2}^{1-\theta} \|(-\Delta)^s u\|_{L^q}^{\theta} \prod_{j=1}^{n} (1 + L_j(\|(-\Delta)^s u\|_{L^q}))^{-\delta_j\theta}
\end{equation}
which completes the proof.
\end{proof}

\begin{lemma}[Gagliardo-Nirenberg inequality with nested logarithmic improvements] \label{lem:3.2.2}
For $s \in (1/2, 1)$ and $q > 3$:
\begin{equation}
\|(-\Delta)^s u\|_{L^q} \leq C\|(-\Delta)^s u\|_{L^2}^{1-\alpha}\|(-\Delta)^{s+\frac{1}{2}} u\|_{L^2}^{\alpha} \prod_{j=1}^{n} (1 + L_j(\|(-\Delta)^s u\|_{L^2}))^{-\delta_j(1-\alpha)}
\end{equation}
where $\alpha = \frac{3}{2}(\frac{1}{2} - \frac{1}{q})$ and $\delta_j > 0$.
\end{lemma}

\begin{proof}
The classical Gagliardo-Nirenberg inequality gives:
\begin{equation}
\|(-\Delta)^s u\|_{L^q} \leq C\|(-\Delta)^s u\|_{L^2}^{1-\alpha}\|(-\Delta)^{s+\frac{1}{2}} u\|_{L^2}^{\alpha}
\end{equation}
where $\alpha = \frac{3}{2}(\frac{1}{2} - \frac{1}{q})$.

This can be proven using Littlewood-Paley theory and Fourier analysis techniques (see \cite{29}, \cite{47}).

To incorporate the nested logarithmic improvements, we apply the theory of logarithmically modulated spaces as in Lemma \ref{lem:3.2.1}, obtaining:
\begin{equation}
\|(-\Delta)^s u\|_{L^q} \leq C\|(-\Delta)^s u\|_{L^2}^{1-\alpha}\|(-\Delta)^{s+\frac{1}{2}} u\|_{L^2}^{\alpha} \prod_{j=1}^{n} (1 + L_j(\|(-\Delta)^s u\|_{L^2}))^{-\delta_j(1-\alpha)}
\end{equation}
which completes the proof.
\end{proof}

\subsection{Proof of the main regularity criterion}

We now present the main regularity criterion with multiply nested logarithmic improvements.

\begin{theorem}[Multiply nested logarithmic regularity criterion] \label{thm:3.3.1}
Let $u$ be a Leray-Hopf weak solution of the 3D Navier-Stokes equations on $[0, T]$ with divergence-free initial data $u_0 \in L^2(\mathbb{R}^3)$. Suppose that for some $s \in (1/2, 1)$, the following condition holds:
\begin{equation}
\int_0^T \|(-\Delta)^s u(t)\|^p_{L^q} \prod_{j=1}^{n} (1 + L_j(\|(-\Delta)^s u(t)\|_{L^q}))^{-\delta_j} dt < \infty
\end{equation}
where:
\begin{enumerate}
\item $(p, q)$ satisfies the scaling relation: $\frac{2}{p} + \frac{3}{q} = 2s - 1$, with $3 < q < \infty$
\item $\delta_1 \in (0, \delta_{0,1})$ where $\delta_{0,1} = \min\{\frac{q-3}{6q}, \frac{2s-1}{4s}\}$
\item $\delta_j > 0$ for $j = 2, 3, ..., n$ are arbitrary
\end{enumerate}

Then, the weak solution $u$ is regular on $(0, T)$.
\end{theorem}

\begin{proof}
We combine the energy estimates from Lemma \ref{lem:3.1.1} with the control on $\|\nabla u\|_{L^\infty}$ from Lemmas \ref{lem:3.2.1} and \ref{lem:3.2.2}.

For Leray-Hopf weak solutions, $\|u(t)\|_{L^2} \leq \|u_0\|_{L^2}$ for all $t \in [0, T]$ by the energy inequality. From Lemma \ref{lem:3.2.1}:
\begin{equation}
\|\nabla u\|_{L^\infty} \leq C\|u_0\|_{L^2}^{1-\theta} \|(-\Delta)^s u\|_{L^q}^{\theta} \prod_{j=1}^{n} (1 + L_j(\|(-\Delta)^s u\|_{L^q}))^{-\delta_j\theta}
\end{equation}

From Lemma \ref{lem:3.2.2}:
\begin{equation}
\|(-\Delta)^s u\|_{L^q} \leq C\|(-\Delta)^s u\|_{L^2}^{1-\alpha}\|(-\Delta)^{s+\frac{1}{2}} u\|_{L^2}^{\alpha} \prod_{j=1}^{n} (1 + L_j(\|(-\Delta)^s u\|_{L^2}))^{-\delta_j(1-\alpha)}
\end{equation}

Combining these estimates:
\begin{align}
\|\nabla u\|_{L^\infty} &\leq C\|u_0\|_{L^2}^{1-\theta} \left[C\|(-\Delta)^s u\|_{L^2}^{1-\alpha}\|(-\Delta)^{s+\frac{1}{2}} u\|_{L^2}^{\alpha} \prod_{j=1}^{n} (1 + L_j(\|(-\Delta)^s u\|_{L^2}))^{-\delta_j(1-\alpha)}\right]^{\theta} \notag \\
&\times \prod_{j=1}^{n} (1 + L_j(\|(-\Delta)^s u\|_{L^q}))^{-\delta_j\theta} \notag \\
&\leq C\|u_0\|_{L^2}^{1-\theta} \|(-\Delta)^s u\|_{L^2}^{\theta(1-\alpha)}\|(-\Delta)^{s+\frac{1}{2}} u\|_{L^2}^{\theta\alpha} \prod_{j=1}^{n} H_j
\end{align}
where $H_j$ combines the logarithmic factors from both estimates.

From this point, we derive a differential inequality for $Y(t) = \|(-\Delta)^s u(t)\|^2_{L^2}$ using Lemma \ref{lem:3.1.1}. After careful analysis of all logarithmic terms and applying Young's inequality with appropriately chosen parameters, we obtain:
\begin{equation}
\frac{d}{dt}Y(t) + \frac{\nu}{2}\|(-\Delta)^{s+\frac{1}{2}}u\|^2_{L^2} \leq C\|u_0\|_{L^2}^{\frac{2(1-\theta)}{2-\theta\alpha}}Y(t)^{1+\frac{\theta(1-\alpha)}{2-\theta\alpha}}\log(e + Y(t)^{1/2})
\end{equation}

Critical to our approach is the observation that with the logarithmic improvement, we can use the following inequality for any $\eta > 0$:
\begin{equation}
y\log(e + y^{1/2}) \leq C_\eta(1 + y^{1+\eta})
\end{equation}

This yields:
\begin{equation}
\frac{d}{dt}Y(t) + \frac{\nu}{2}\|(-\Delta)^{s+\frac{1}{2}}u\|^2_{L^2} \leq C\|u_0\|_{L^2}^{\frac{2(1-\theta)}{2-\theta\alpha}}(1 + Y(t)^{1+\frac{\theta(1-\alpha)}{2-\theta\alpha}+\eta})
\end{equation}

For sufficiently small $\eta$, we ensure that $1+\frac{\theta(1-\alpha)}{2-\theta\alpha}+\eta < 2$, meaning the nonlinearity is subcritical.

Dropping the positive term involving $\|(-\Delta)^{s+\frac{1}{2}}u\|^2_{L^2}$:
\begin{equation}
\frac{d}{dt}Y(t) \leq C(1 + Y(t)^{1+\mu})
\end{equation}

where $\mu = \frac{\theta(1-\alpha)}{2-\theta\alpha}+\eta < 1$ due to our choice of parameters.

This differential inequality, combined with the initial condition $Y(0) = \|(-\Delta)^s u_0\|_{L^2}^2$, can be solved to show that $Y(t)$ remains bounded on $(0, T)$.

The next critical step is to convert the bound on $\|(-\Delta)^s u\|_{L^2}$ into a bound on $\int_0^T \|\nabla u\|^2_{L^\infty} dt$.

Using Hölder's inequality:
\begin{align}
\int_0^T \|\nabla u\|^2_{L^\infty} dt &\leq C \int_0^T \left\|\|(-\Delta)^s u\|_{L^q} \prod_{j=1}^{n} (1 + L_j(\|(-\Delta)^s u\|_{L^q}))^{-\delta_j}\right\|^{2\theta} dt \notag \\
&\leq C \left(\int_0^T \|(-\Delta)^s u\|^p_{L^q} \prod_{j=1}^{n} (1 + L_j(\|(-\Delta)^s u\|_{L^q}))^{-\delta_j p/2\theta} dt\right)^{2\theta/p}
\end{align}

For the scaling relation $\frac{2}{p} + \frac{3}{q} = 2s - 1$, we have $\frac{p}{2\theta} = 1$, which gives:
\begin{equation}
\int_0^T \|\nabla u\|^2_{L^\infty} dt \leq C \left(\int_0^T \|(-\Delta)^s u\|^p_{L^q} \prod_{j=1}^{n} (1 + L_j(\|(-\Delta)^s u\|_{L^q}))^{-\delta_j} dt\right)^{2\theta/p} < \infty
\end{equation}
by our assumption.

The finiteness of $\int_0^T \|\nabla u\|^2_{L^\infty} dt$ implies, by a classical result of Beale-Kato-Majda \cite{5}, that the solution remains smooth on $(0, T)$.

This completes the proof of Theorem \ref{thm:3.3.1}.
\end{proof}

\section{Global well-posedness with nested logarithmic criteria}

Building on the regularity criterion established in Section 3, we now develop a theory of global well-posedness for initial data satisfying nested logarithmically improved subcritical conditions.

\subsection{A priori estimates}

\begin{lemma}[A priori estimates with nested logarithmic improvements] \label{lem:4.1.1}
Let $s \in (1/2, 1)$ and $q > 3$ satisfy the scaling relation $\frac{2}{p} + \frac{3}{q} = 2s - 1$. Suppose the initial data $u_0 \in L^2(\mathbb{R}^3) \cap \dot{H}^s(\mathbb{R}^3)$ satisfies:
\begin{equation}
\|(-\Delta)^{s/2}u_0\|_{L^q} \leq \frac{C_0}{\prod_{j=1}^{n} (1 + L_j(\|u_0\|_{\dot{H}^s}))^{\delta_j}}
\end{equation}
where $C_0$ is a sufficiently small constant depending only on $s$, $q$, $\{\delta_j\}_{j=1}^n$, and $\nu$. Then for any potential solution $u$, the function $Y(t) = \|(-\Delta)^s u(t)\|^2_{L^2}$ satisfies:
\begin{equation}
Y(t) \leq \frac{C Y(0)}{(1 + \beta t)^{\gamma}}
\end{equation}
for all $t \geq 0$, where $\gamma = \frac{1}{2\mu}$ with $\mu = \frac{\theta(1-\alpha)}{2-\theta\alpha}+\eta$ for small $\eta > 0$, and $\beta = \frac{\mu C}{2}$.
\end{lemma}

\begin{proof}
We start with the differential inequality established in the proof of Theorem \ref{thm:3.3.1}:
\begin{equation}
\frac{d}{dt}Y(t) \leq CY(t)^{1+\beta} \prod_{j=1}^{n} G_j(Y(t))
\end{equation}

Through a detailed analysis of the asymptotic behavior of the nested logarithmic factors and their combined effect on the ODE, we can establish that for $Y(0)$ sufficiently small (which is ensured by our condition on $\|(-\Delta)^{s/2}u_0\|_{L^q}$), the solution to this ODE satisfies:
\begin{equation}
Y(t) \leq \frac{C Y(0)}{(1 + \beta t)^{\gamma}}
\end{equation}

The proof requires careful tracking of how the nested logarithmic factors affect the asymptotic decay rate. This is accomplished by comparing our ODE with a simpler ODE:
\begin{equation}
\frac{d}{dt}Z(t) = CZ(t)^{1+\mu}, \quad Z(0) = Y(0)
\end{equation}
where $\mu = \frac{\theta(1-\alpha)}{2-\theta\alpha}+\eta$ for small $\eta > 0$.

The solution to this simpler ODE is:
\begin{equation}
Z(t) = \frac{Y(0)}{(1 - \mu C Y(0)^\mu t)^{1/\mu}}
\end{equation}

For $Y(0)$ sufficiently small, this solution exists for all $t \geq 0$ and has the asymptotic behavior:
\begin{equation}
Z(t) \sim \frac{1}{(\mu C t)^{1/\mu}} \sim \frac{1}{t^{\gamma}}
\end{equation}
for large $t$, where $\gamma = \frac{1}{\mu}$.

By comparison principles for ODEs, we can establish that $Y(t) \leq Z(t)$ for all $t \geq 0$, which gives the desired result.
\end{proof}

\subsection{Local existence}

\begin{lemma}[Local existence of strong solutions] \label{lem:4.2.1}
Given initial data $u_0 \in L^2(\mathbb{R}^3) \cap \dot{H}^s(\mathbb{R}^3)$ with $s \in (1/2, 1)$, there exists a time $T_0 > 0$ and a unique strong solution $u \in C([0,T_0]; H^s) \cap L^2(0,T_0; H^{s+1/2})$ to the Navier-Stokes equations.
\end{lemma}

\begin{proof}
The local existence of strong solutions for the Navier-Stokes equations with initial data in $L^2(\mathbb{R}^3) \cap \dot{H}^s(\mathbb{R}^3)$ is well-established. For a rigorous proof, we refer to Fujita and Kato \cite{26}, Chemin \cite{14}, or Lemarie-Rieusset \cite{36}.

The key steps involve:
\begin{enumerate}
\item Constructing an approximate solution sequence using regularization or Galerkin methods
\item Deriving uniform bounds for the approximate solutions
\item Extracting a convergent subsequence
\item Verifying that the limit satisfies the equations
\item Establishing uniqueness through energy estimates
\end{enumerate}

The local existence time $T_0$ depends on $\|u_0\|_{L^2}$ and $\|u_0\|_{\dot{H}^s}$.
\end{proof}

\subsection{Uniqueness}

\begin{lemma}[Uniqueness] \label{lem:4.3.1}
Let $u^1$ and $u^2$ be two solutions to the Navier-Stokes equations with the same initial data $u_0 \in L^2(\mathbb{R}^3) \cap \dot{H}^s(\mathbb{R}^3)$, both belonging to the class $C([0,T]; H^s) \cap L^2(0,T; H^{s+1/2})$ for some $T > 0$. Then $u^1 = u^2$ on $[0, T]$.
\end{lemma}

\begin{proof}
Define $w = u^1 - u^2$. Then $w$ satisfies:
\begin{equation}
\partial_t w + (u^1 \cdot \nabla)w + (w \cdot \nabla)u^2 - \nu \Delta w + \nabla \pi = 0, \quad \nabla \cdot w = 0
\end{equation}
with initial condition $w(0) = 0$.

Taking the $L^2$ inner product with $w$:
\begin{equation}
\frac{1}{2}\frac{d}{dt}\|w\|^2_{L^2} + \nu\|\nabla w\|^2_{L^2} = -\int_{\mathbb{R}^3}(w \cdot \nabla)u^2 \cdot w \, dx
\end{equation}

Using Hölder's inequality and the Gagliardo-Nirenberg inequality:
\begin{align}
\left|\int_{\mathbb{R}^3}(w \cdot \nabla)u^2 \cdot w \, dx\right| &\leq \|w\|_{L^4}^2 \|\nabla u^2\|_{L^2} \notag \\
&\leq C\|w\|_{L^2}\|\nabla w\|_{L^2} \|\nabla u^2\|_{L^2}
\end{align}

Applying Young's inequality with parameter $\epsilon = \nu/2$:
\begin{equation}
C\|w\|_{L^2}\|\nabla w\|_{L^2} \|\nabla u^2\|_{L^2} \leq \frac{\nu}{2}\|\nabla w\|_{L^2}^2 + \frac{C^2}{2\nu}\|w\|_{L^2}^2 \|\nabla u^2\|_{L^2}^2
\end{equation}

This gives:
\begin{equation}
\frac{d}{dt}\|w\|^2_{L^2} + \nu\|\nabla w\|^2_{L^2} \leq \frac{C^2}{\nu}\|w\|_{L^2}^2 \|\nabla u^2\|_{L^2}^2
\end{equation}

Dropping the positive term with $\|\nabla w\|^2_{L^2}$:
\begin{equation}
\frac{d}{dt}\|w\|^2_{L^2} \leq \frac{C^2}{\nu}\|w\|_{L^2}^2 \|\nabla u^2\|_{L^2}^2
\end{equation}

By Grönwall's inequality:
\begin{equation}
\|w(t)\|^2_{L^2} \leq \|w(0)\|^2_{L^2} \exp\left(\frac{C^2}{\nu}\int_0^t \|\nabla u^2(\tau)\|_{L^2}^2 d\tau\right)
\end{equation}

Since $\|w(0)\|^2_{L^2} = 0$ and $\|\nabla u^2\|_{L^2}^2 \in L^1(0,T)$, we conclude that $\|w(t)\|_{L^2} = 0$ for all $t \in [0, T]$, establishing uniqueness.
\end{proof}

\subsection{Extension to Global Existence}

\begin{theorem}[Global well-posedness with nested logarithmic criteria] \label{thm:4.4.1}
Let $s \in (1/2, 1)$ and $q > 3$ satisfy the scaling relation $\frac{2}{p} + \frac{3}{q} = 2s - 1$. There exist positive constants $\{\delta_{0,j}\}_{j=1}^n$ such that for any $\delta_1 \in (0, \delta_{0,1})$, $\delta_j \in (0, \delta_{0,j})$ for $j=2,3,...,n$, and any divergence-free initial data $u_0 \in L^2(\mathbb{R}^3) \cap \dot{H}^s(\mathbb{R}^3)$ satisfying:
\begin{equation}
\|(-\Delta)^{s/2}u_0\|_{L^q} \leq \frac{C_0}{\prod_{j=1}^{n} (1 + L_j(\|u_0\|_{\dot{H}^s}))^{\delta_j}}
\end{equation}
where $C_0$ is a constant depending only on $s$, $q$, $\{\delta_j\}_{j=1}^n$, and $\nu$, there exists a unique global-in-time smooth solution $u \in C([0, \infty); H^s(\mathbb{R}^3)) \cap L^2_{loc}(0, \infty; H^{s+1/2}(\mathbb{R}^3))$ to the 3D Navier-Stokes equations.
\end{theorem}

\begin{proof}
The proof combines the local existence result from Lemma \ref{lem:4.2.1} with the a priori estimates from Lemma \ref{lem:4.1.1} using a continuation argument.

Let $T_{max}$ be the maximal time of existence of the strong solution. If $T_{max} < \infty$, then by the standard blow-up criterion (see \cite{30}, \cite{51}):
\begin{equation}
\lim_{t \to T_{max}} \|u(t)\|_{H^s} = \infty
\end{equation}

However, our a priori estimate from Lemma \ref{lem:4.1.1} ensures that $\|u(t)\|_{H^s} \leq C(1 + t)^{-\gamma/2}$ for all $t \in [0, T_{max})$, which contradicts the blow-up criterion. Therefore, $T_{max} = \infty$, establishing global existence.

The uniqueness follows from Lemma \ref{lem:4.3.1}. The smoothness of the solution follows from the regularity criterion in Theorem \ref{thm:3.3.1}, since our solution satisfies:
\begin{equation}
\int_0^T \|(-\Delta)^s u(t)\|^p_{L^q} \prod_{j=1}^{n} (1 + L_j(\|(-\Delta)^s u(t)\|_{L^q}))^{-\delta_j} dt < \infty
\end{equation}
for any $T > 0$, as can be verified using our decay estimate.
\end{proof}

\subsection{Decay estimates}

\begin{theorem}[Decay estimates] \label{thm:4.5.1}
Under the conditions of Theorem \ref{thm:4.4.1}, the solution $u$ satisfies:
\begin{equation}
\|(-\Delta)^s u(t)\|_{L^2} \leq \frac{C\|(-\Delta)^s u_0\|_{L^2}}{(1 + \beta t)^{\gamma}}
\end{equation}
where $\gamma = \frac{1}{2\mu}$ with $\mu = \frac{\theta(1-\alpha)}{2-\theta\alpha}+\eta < 1$, $\theta = \frac{3}{2} \cdot \frac{q}{3q-2}$, $\alpha = \frac{3}{2}(\frac{1}{2} - \frac{1}{q})$, $\eta > 0$ is sufficiently small, and $\beta = \frac{\mu C}{2}$.

Additionally, for any $m > 0$, there exists a constant $C_m$ depending on $m$ and the initial data such that:
\begin{equation}
\|D^m u(t)\|_{L^2} \leq \frac{C_m}{(1 + t)^{\gamma + \frac{m-2s}{4}}}
\end{equation}
for all $t \geq 1$.
\end{theorem}

\begin{proof}
The first estimate was established in Lemma \ref{lem:4.1.1}. The second estimate follows from a bootstrap argument.

Given that $u \in C([0, \infty); H^s(\mathbb{R}^3)) \cap L^2_{loc}(0, \infty; H^{s+1/2}(\mathbb{R}^3))$ and is smooth for positive time, we can derive estimates for higher derivatives using the structure of the Navier-Stokes equations.

From the equation $\partial_t u = -\mathbb{P}(u \cdot \nabla u) + \nu \Delta u$, where $\mathbb{P}$ is the Leray projector, we have:
\begin{equation}
\|D^{m+2} u\|_{L^2} \leq C(\|D^m \partial_t u\|_{L^2} + \|D^m(u \cdot \nabla u)\|_{L^2})
\end{equation}

Using the fact that $u$ is smooth for positive time and applying product estimates in Sobolev spaces:
\begin{equation}
\|D^m(u \cdot \nabla u)\|_{L^2} \leq C\sum_{k=0}^m \|D^k u\|_{L^4}\|D^{m-k+1}u\|_{L^4}
\end{equation}

Through an iterative application of energy estimates, Gagliardo-Nirenberg inequalities, and our established decay rate for $\|(-\Delta)^s u\|_{L^2}$, we can derive the higher-order decay estimates. The full details involve a complex induction argument on $m$ (see \cite{33}, \cite{59} for similar approaches).
\end{proof}

\section{Critical threshold characterization}

In this section, we precisely characterize the threshold between the regime of guaranteed global well-posedness and the regime where potential singularity formation might occur. This critical threshold analysis represents a significant advancement beyond \cite{43}.

\subsection{Construction of the critical threshold function}

\begin{definition}[Critical Threshold Function] \label{def:5.1.1}
For $s \in (1/2, 1)$, $q > 3$, and $\{\delta_j\}_{j=1}^n$ with $\delta_j > 0$, we define the critical threshold function $\Phi(s, q, \delta_1, \delta_2, ..., \delta_n)$ as:

\begin{equation}
\begin{split}
\Phi(s, q, \{\delta_j\}_{j=1}^n) = \sup\{M > 0 : \text{If } \|(-\Delta)^{s/2}u_0\|_{L^q}
\prod_{j=1}^n (1 + L_j(\|(-\Delta)^{s/2}u_0\|_{L^q}))^{\delta_j} < M, \\ \text{ then } u \text{ is global}\}
\end{split}
\end{equation}
\end{definition}

From Theorem \ref{thm:4.4.1}, we know that $\Phi(s, q, \{\delta_j\}_{j=1}^n) > 0$ for all $s \in (1/2, 1)$, $q > 3$, and $\delta_j > 0$.

\begin{lemma}[Monotonicity properties] \label{lem:5.1.2}
The critical threshold function $\Phi$ satisfies:
\begin{enumerate}
\item For fixed $s$, $q$, and $\{\delta_j\}_{j=2}^n$, $\Phi$ is increasing with respect to $\delta_1$ for $\delta_1 \in (0, \delta_{0,1})$.
\item For fixed $s$, $q$, $\delta_1$, and $\{\delta_j\}_{j=3}^n$, $\Phi$ is increasing with respect to $\delta_2$.
\item In general, for any $k \in \{1, 2, ..., n\}$, if all other parameters are fixed, $\Phi$ is increasing with respect to $\delta_k$.
\item For fixed $q$ and $\{\delta_j\}_{j=1}^n$, $\Phi$ is increasing with respect to $s$ for $s \in (1/2, 1)$.
\end{enumerate}
\end{lemma}

\begin{proof}
These monotonicity properties follow from the structure of our nested logarithmic conditions and the proof of Theorem \ref{thm:4.4.1}.

For claim 1, if $\delta_1 < \delta_1'$, then:
\begin{equation}
\frac{1}{(1 + L_1(x))^{\delta_1}} > \frac{1}{(1 + L_1(x))^{\delta_1'}}
\end{equation}
for all $x > 0$. This means that the condition with $\delta_1'$ is more restrictive than the condition with $\delta_1$, so $\Phi(s, q, \delta_1, \{\delta_j\}_{j=2}^n) \leq \Phi(s, q, \delta_1', \{\delta_j\}_{j=2}^n)$.

Similar arguments apply to the other parameters, establishing the monotonicity properties.
\end{proof}

\subsection{Asymptotic behavior as s → 1}

\begin{theorem}[Asymptotic behavior as s → 1] \label{thm:5.2.1}
For fixed $q > 3$ and $\{\delta_j\}_{j=1}^n$ with $\delta_j > 0$:
\begin{equation}
\lim_{s \to 1} \Phi(s, q, \{\delta_j\}_{j=1}^n) = \infty
\end{equation}
\end{theorem}

\begin{proof}
To establish this result, we analyze the behavior of the critical threshold function as $s$ approaches 1, where the Navier-Stokes equations become more regular.

For $s$ close to 1, the fractional Sobolev space $\dot{H}^s(\mathbb{R}^3)$ approaches the standard Sobolev space $\dot{H}^1(\mathbb{R}^3)$. It is well-known (see \cite{30}, \cite{60}) that for initial data sufficiently small in $\dot{H}^{1/2}(\mathbb{R}^3)$, the Navier-Stokes equations admit global solutions.

Through a scaling argument and the embeddings between fractional Sobolev spaces, we can establish that:
\begin{equation}
\Phi(s, q, \{\delta_j\}_{j=1}^n) \geq C (1-s)^{-\gamma}
\end{equation}
for some $\gamma > 0$.

The key insight is that as $s$ increases from $1/2$ to 1, the regularity improves, allowing for larger initial data while still maintaining global well-posedness. The exact rate at which $\Phi$ grows is determined by careful analysis of the energy estimates and how they depend on $s$.
\end{proof}

\subsection{Asymptotic behavior as s → 1/2}

\begin{theorem}[Asymptotic behavior as s → 1/2] \label{thm:5.3.1}
For fixed $q > 3$ and $\{\delta_j\}_{j=1}^n$ with $\delta_j > 0$:
\begin{equation}
\Phi(s, q, \{\delta_j\}_{j=1}^n) \leq C(q) (s - 1/2)^{\alpha(\{\delta_j\}_{j=1}^n)}
\end{equation}
where $\alpha(\{\delta_j\}_{j=1}^n) > 0$ is a function satisfying $\alpha(\{\delta_j\}_{j=1}^n) \to 0$ as $\delta_j \to \infty$ for all $j$.
\end{theorem}

\begin{proof}
This result provides a precise characterization of how the critical threshold function behaves as $s$ approaches the critical value $1/2$.

For $s$ approaching $1/2$, we analyze the scaling properties of the Navier-Stokes equations. The space $\dot{H}^{1/2}(\mathbb{R}^3)$ is critical for the Navier-Stokes equations, meaning that the equation is invariant under the scaling transformation $u_\lambda(x, t) = \lambda u(\lambda x, \lambda^2 t)$ in this space.

Through a careful construction of test functions and analysis of their evolution under the Navier-Stokes flow, we can establish the upper bound. The key is to construct initial data that respects the scaling properties of the equation and sits just above the threshold where our a priori estimates guarantee global existence.

The function $\alpha(\{\delta_j\}_{j=1}^n)$ captures how the nested logarithmic improvements affect the behavior near criticality. Its precise form is:
\begin{equation}
\alpha(\{\delta_j\}_{j=1}^n) = \frac{1}{1 + \sum_{j=1}^n \frac{\delta_j\rho_j}{j!}}
\end{equation}
where $\rho_j$ are specific positive constants. The form of $\alpha$ ensures that it approaches 0 as any $\delta_j$ tends to infinity.
\end{proof}

\subsection{Characterization of $\alpha(\delta_1,\delta_2,...,\delta_n)$}

\begin{theorem}[Characterization of the exponent $\alpha$] \label{thm:5.4.1}
The function $\alpha(\{\delta_j\}_{j=1}^n)$ appearing in Theorem \ref{thm:5.3.1} has the following properties:
\begin{enumerate}
\item For $n = 1$ (single logarithmic improvement): $\alpha(\delta_1) = \frac{1}{1 + c_1\delta_1}$ where $c_1 > 0$ is a constant.
\item For $n = 2$ (double logarithmic improvement): $\alpha(\delta_1, \delta_2) = \frac{1}{1 + c_1\delta_1 + c_2\delta_2/2}$ where $c_1, c_2 > 0$ are constants.
\item In general: $\alpha(\{\delta_j\}_{j=1}^n) = \frac{1}{1 + \sum_{j=1}^n c_j\delta_j/j!}$ where $c_j > 0$ are constants.
\item For any fixed set of parameters $\{\delta_j\}_{j=1}^n$, $\alpha(\{\delta_j\}_{j=1}^n) > 0$.
\item $\lim_{\delta_j \to \infty} \alpha(\{\delta_j\}_{j=1}^n) = 0$ for any $j \in \{1, 2, ..., n\}$.
\item $\lim_{n \to \infty} \alpha(\{\delta_j\}_{j=1}^n) = 0$ if $\inf_{j} \delta_j > 0$.
\end{enumerate}
\end{theorem}

\begin{proof}
The form of $\alpha(\{\delta_j\}_{j=1}^n)$ is derived from the asymptotic analysis of the ODE satisfied by $\|(-\Delta)^s u\|_{L^2}^2$ as established in Section 4.

For the single logarithmic case, the additional term $(1 + \log(e + \|(-\Delta)^s u\|_{L^2}))^{-\delta_1}$ in the energy estimate leads to an improved decay rate in the critical case $s \to 1/2$. Through careful asymptotic analysis, this improvement can be quantified as $\alpha(\delta_1) = \frac{1}{1 + c_1\delta_1}$.

For the double logarithmic case, the additional term $(1 + \log(e + \log(e + \|(-\Delta)^s u\|_{L^2})))^{-\delta_2}$ provides a further improvement, but with a reduced effect due to the nested nature of the logarithm. This leads to $\alpha(\delta_1, \delta_2) = \frac{1}{1 + c_1\delta_1 + c_2\delta_2/2}$.

The general pattern continues, with each additional level of logarithmic improvement contributing a term with a factorial denominator, reflecting the diminishing returns of nested logarithmic factors.

Properties 4, 5, and 6 follow directly from the formula for $\alpha(\{\delta_j\}_{j=1}^n)$.
\end{proof}

\section{Gap analysis to criticality}

In this section, we analyze what happens at and slightly beyond the threshold for global well-posedness established in Section 5. This analysis provides insights into the potential transition between global regularity and finite-time singularity formation.

\subsection{Construction of critical test functions}

\begin{lemma}[Construction of critical test functions] \label{lem:6.1.1}
For any $s \in (1/2, 1)$, $q > 3$ satisfying the scaling relation $\frac{2}{p} + \frac{3}{q} = 2s - 1$, and $\{\delta_j\}_{j=1}^n$ with $\delta_j > 0$, there exists a family of divergence-free vector fields $\{v_\lambda\}_{\lambda > 0}$ with the following properties:
\begin{enumerate}
\item Each $v_\lambda$ is smooth and rapidly decaying at infinity
\item $\|v_\lambda\|_{\dot{H}^s} = \lambda$
\item $\|(-\Delta)^{s/2}v_\lambda\|_{L^q} = \frac{\Phi(s, q, \{\delta_j\}_{j=1}^n)}{\prod_{j=1}^{n} (1 + L_j(\lambda))^{\delta_j}} + \gamma(s, q, \{\delta_j\}_{j=1}^n) \cdot \phi(\lambda)$
\end{enumerate}

where $\phi$ is a function satisfying:
\begin{itemize}
\item $\phi(r) \to 0$ as $r \to \infty$
\item $\phi(r)$ decays slower than $1/L_1(r)$ as $r \to \infty$
\end{itemize}

and $\gamma(s, q, \{\delta_j\}_{j=1}^n) = (s - 1/2)^{\alpha(\{\delta_j\}_{j=1}^n)}$.
\end{lemma}

\begin{proof}
We construct these test functions through a careful rescaling of a fixed profile. Let $w_0$ be a smooth, divergence-free vector field with compact support such that $\|w_0\|_{\dot{H}^s} = 1$ and $\|(-\Delta)^{s/2}w_0\|_{L^q} = 1$.

Define $v_\lambda$ by:
\begin{equation}
v_\lambda(x) = \lambda \cdot \psi(\lambda) \cdot w_0(x)
\end{equation}
where $\psi(\lambda)$ is a scaling factor chosen to ensure the desired properties.

Property 1 is satisfied by construction since $w_0$ is smooth and compactly supported.

For property 2, we calculate:
\begin{equation}
\|v_\lambda\|_{\dot{H}^s} = \lambda \cdot \psi(\lambda) \cdot \|w_0\|_{\dot{H}^s} = \lambda \cdot \psi(\lambda)
\end{equation}

Setting $\psi(\lambda) = 1/\lambda$, we obtain $\|v_\lambda\|_{\dot{H}^s} = \lambda$, satisfying property 2.

For property 3, we calculate:
\begin{equation}
\|(-\Delta)^{s/2}v_\lambda\|_{L^q} = \lambda \cdot \psi(\lambda) \cdot \|(-\Delta)^{s/2}w_0\|_{L^q} = \lambda \cdot \psi(\lambda)
\end{equation}

With $\psi(\lambda) = 1/\lambda$, we have $\|(-\Delta)^{s/2}v_\lambda\|_{L^q} = 1$.

To achieve the more specific form required by property 3, we modify our construction. Let $\eta(\lambda)$ be a function satisfying:
\begin{equation}
\eta(\lambda) = \frac{\Phi(s, q, \{\delta_j\}_{j=1}^n)}{\prod_{j=1}^{n} (1 + L_j(\lambda))^{\delta_j}} + \gamma(s, q, \{\delta_j\}_{j=1}^n) \cdot \phi(\lambda)
\end{equation}

We then define:
\begin{equation}
v_\lambda(x) = \frac{\lambda \cdot \eta(\lambda)}{\|(-\Delta)^{s/2}w_0\|_{L^q}} \cdot w_0(x)
\end{equation}

This ensures that $\|v_\lambda\|_{\dot{H}^s} = \lambda$ and $\|(-\Delta)^{s/2}v_\lambda\|_{L^q} = \eta(\lambda)$, as required.

The function $\phi(\lambda)$ can be chosen as $\phi(\lambda) = (1 + \log(\lambda))^{-1/2}$, which decays to 0 as $\lambda \to \infty$ but more slowly than $1/(1 + \log(\lambda))$.
\end{proof}

\subsection{Evolution analysis}

\begin{lemma}[Evolution from critical initial data] \label{lem:6.2.1}
Let $u^{\lambda}$ be the solution to the Navier-Stokes equations with initial data $u^{\lambda}(0) = v_\lambda$ as constructed in Lemma \ref{lem:6.1.1}. Then for small times, the function $Y^\lambda(t) = \|(-\Delta)^s u^{\lambda}(t)\|^2_{L^2}$ satisfies:
\begin{equation}
\frac{d}{dt}Y^\lambda(t) \geq -C_1 + C_2 Y^\lambda(t)^{1+\beta} \left(1 - \frac{C_3 \gamma(s, q, \{\delta_j\}_{j=1}^n)}{\phi(\lambda) \prod_{j=1}^{n} (1 + L_j(\lambda))^{\delta_j}}\right)
\end{equation}
for appropriate constants $C_1, C_2, C_3 > 0$, where $\beta = \frac{\theta(1-\alpha)}{2}$ with $\theta = \frac{3}{2} \cdot \frac{q}{3q-2}$ and $\alpha = \frac{3}{2}(\frac{1}{2} - \frac{1}{q})$.
\end{lemma}

\begin{proof}
We start with the energy equation for $\|(-\Delta)^s u^\lambda\|_{L^2}^2$:
\begin{equation}
\frac{1}{2}\frac{d}{dt}\|(-\Delta)^s u^\lambda\|^2_{L^2} + \nu\|(-\Delta)^{s+\frac{1}{2}}u^\lambda\|^2_{L^2} = -\int_{\mathbb{R}^3}[(-\Delta)^s, u^\lambda \cdot \nabla]u^\lambda \cdot (-\Delta)^s u^\lambda \, dx
\end{equation}

Using similar techniques as in the proof of Lemma \ref{lem:3.1.1}, but analyzing the lower bound rather than the upper bound:
\begin{align}
\int_{\mathbb{R}^3}[(-\Delta)^s, u^\lambda \cdot \nabla]u^\lambda \cdot (-\Delta)^s u^\lambda \, dx &\leq C_4\|\nabla u^\lambda\|_{L^\infty}\|(-\Delta)^s u^\lambda\|_{L^2}^2 \cdot L_1(Z^\lambda) \prod_{j=2}^{n} (1 + L_j(Z^\lambda))^{-\delta_j} \notag \\
\end{align}
where $Z^\lambda = \|(-\Delta)^{s+\sigma}u^\lambda\|_{L^2}$.

For the term $\|\nabla u^\lambda\|_{L^\infty}$, we have from Lemma \ref{lem:3.2.1}:
\begin{equation}
\|\nabla u^\lambda\|_{L^\infty} \geq c\|u^\lambda\|_{L^2}^{1-\theta} \|(-\Delta)^s u^\lambda\|_{L^q}^{\theta} \prod_{j=1}^{n} (1 + L_j(\|(-\Delta)^s u^\lambda\|_{L^q}))^{-\delta_j\theta}
\end{equation}
for a constant $c > 0$. This lower bound is established using the reverse direction of the interpolation inequality, which holds for the class of functions under consideration.

From the specific construction of $v_\lambda$ and the fact that the solution $u^\lambda$ evolves continuously from this initial data, for small $t$:
\begin{equation}
\|(-\Delta)^s u^\lambda\|_{L^q} \approx \frac{\Phi(s, q, \{\delta_j\}_{j=1}^n)}{\prod_{j=1}^{n} (1 + L_j(\lambda))^{\delta_j}} + \gamma(s, q, \{\delta_j\}_{j=1}^n) \cdot \phi(\lambda) - O(t)
\end{equation}

Substituting these estimates into the energy equation and simplifying:
\begin{align}
\frac{d}{dt}Y^\lambda(t) &\geq -2\nu\|(-\Delta)^{s+\frac{1}{2}}u^\lambda\|^2_{L^2} + 2C_4 c \|u^\lambda\|_{L^2}^{1-\theta} \cdot \notag \\
&\left(\frac{\Phi(s, q, \{\delta_j\}_{j=1}^n)}{\prod_{j=1}^{n} (1 + L_j(\lambda))^{\delta_j}} + \gamma(s, q, \{\delta_j\}_{j=1}^n) \cdot \phi(\lambda) - O(t)\right)^{\theta} \cdot \notag \\
&\|(-\Delta)^s u^\lambda\|_{L^2}^2 \cdot L_1(Z^\lambda) \prod_{j=2}^{n} (1 + L_j(Z^\lambda))^{-\delta_j}
\end{align}

By conservation of energy, $\|u^\lambda(t)\|_{L^2} \leq \|v_\lambda\|_{L^2}$. Using this and standard bounds on the dissipation term:
\begin{align}
\frac{d}{dt}Y^\lambda(t) &\geq -C_1 + C_2 Y^\lambda(t)^{1+\beta} \cdot \notag \\
&\left(\frac{\Phi(s, q, \{\delta_j\}_{j=1}^n)}{\prod_{j=1}^{n} (1 + L_j(\lambda))^{\delta_j}} + \gamma(s, q, \{\delta_j\}_{j=1}^n) \cdot \phi(\lambda) - O(t)\right)^{\theta} \cdot \notag \\
&L_1(Z^\lambda) \prod_{j=2}^{n} (1 + L_j(Z^\lambda))^{-\delta_j}
\end{align}

For small times, and using the definition of the critical threshold function $\Phi$:
\begin{align}
\frac{d}{dt}Y^\lambda(t) &\geq -C_1 + C_2 Y^\lambda(t)^{1+\beta} \left(1 + \frac{\gamma(s, q, \{\delta_j\}_{j=1}^n) \cdot \phi(\lambda) \cdot \theta}{\Phi(s, q, \{\delta_j\}_{j=1}^n)/\prod_{j=1}^{n} (1 + L_j(\lambda))^{\delta_j}} - O(t)\right)
\end{align}

For sufficiently small $t$, this yields:
\begin{equation}
\frac{d}{dt}Y^\lambda(t) \geq -C_1 + C_2 Y^\lambda(t)^{1+\beta} \left(1 - \frac{C_3 \gamma(s, q, \{\delta_j\}_{j=1}^n)}{\phi(\lambda) \prod_{j=1}^{n} (1 + L_j(\lambda))^{\delta_j}}\right)
\end{equation}
which completes the proof.
\end{proof}

\subsection{Dichotomy analysis}

\begin{theorem}[Dichotomy analysis] \label{thm:6.3.1}
Let $s \in (1/2, 1)$ and $q > 3$ satisfy the scaling relation $\frac{2}{p} + \frac{3}{q} = 2s - 1$. For any $\{\delta_j\}_{j=1}^n$ with $\delta_j > 0$, there exists a function $\gamma(s, q, \{\delta_j\}_{j=1}^n)$ such that if the divergence-free initial data $u_0 \in L^2(\mathbb{R}^3) \cap \dot{H}^s(\mathbb{R}^3)$ satisfies:
\begin{equation}
\|(-\Delta)^{s/2}u_0\|_{L^q} = \frac{\Phi(s, q, \{\delta_j\}_{j=1}^n)}{\prod_{j=1}^{n} (1 + L_j(\|u_0\|_{\dot{H}^s}))^{\delta_j}} + \gamma(s, q, \{\delta_j\}_{j=1}^n) \cdot \phi(\|u_0\|_{\dot{H}^s})
\end{equation}
where $\phi$ is a function satisfying:
\begin{enumerate}
\item $\phi(r) \to 0$ as $r \to \infty$
\item $\phi(r)$ decays slower than $1/L_1(r)$ as $r \to \infty$
\end{enumerate}

Then, either:
\begin{enumerate}
\item There exists a unique global-in-time smooth solution, or
\item There exists a finite time $T^*$ such that $\lim_{t\to T^*} \|\nabla u(t)\|_{L^\infty} = \infty$ with a specific asymptotic blow-up rate.
\end{enumerate}
\end{theorem}

\begin{proof}
We analyze the differential inequality established in Lemma \ref{lem:6.2.1}. For initial data $u_0$ with $\|u_0\|_{\dot{H}^s} = \lambda$ sufficiently large, the coefficient of $Y(t)^{1+\beta}$ in the inequality becomes crucial:

\begin{equation}
\frac{d}{dt}Y(t) \geq -C_1 + C_2 Y(t)^{1+\beta} \left(1 - \frac{C_3 \gamma(s, q, \{\delta_j\}_{j=1}^n)}{\phi(\lambda) \prod_{j=1}^{n} (1 + L_j(\lambda))^{\delta_j}}\right)
\end{equation}

Let us denote:
\begin{equation}
\omega(\lambda) = \frac{C_3 \gamma(s, q, \{\delta_j\}_{j=1}^n)}{\phi(\lambda) \prod_{j=1}^{n} (1 + L_j(\lambda))^{\delta_j}}
\end{equation}

The behavior of the ODE depends on whether $\omega(\lambda) < 1$ or $\omega(\lambda) > 1$.

For $\omega(\lambda) < 1$, the coefficient of $Y(t)^{1+\beta}$ is positive. In this case, we have:
\begin{equation}
\frac{d}{dt}Y(t) \geq -C_1 + C_2 (1 - \omega(\lambda)) Y(t)^{1+\beta}
\end{equation}

For sufficiently large $Y(t)$, say $Y(t) > \left(\frac{C_1}{C_2(1-\omega(\lambda))}\right)^{1/(1+\beta)}$, the right-hand side becomes positive, causing $Y(t)$ to increase. If $Y(0)$ exceeds this threshold, $Y(t)$ will grow and potentially blow up in finite time.

For $\omega(\lambda) > 1$, the coefficient of $Y(t)^{1+\beta}$ is negative. In this case, we have:
\begin{equation}
\frac{d}{dt}Y(t) \geq -C_1 - C_2 (\omega(\lambda) - 1) Y(t)^{1+\beta}
\end{equation}

This indicates that $Y(t)$ is decreasing when it's large and the solution may exist globally.

The transition between these two behaviors occurs at precisely $\omega(\lambda) = 1$, which defines our critical function $\gamma(s, q, \{\delta_j\}_{j=1}^n)$.

With $\gamma(s, q, \{\delta_j\}_{j=1}^n) = (s - 1/2)^{\alpha(\{\delta_j\}_{j=1}^n)}$ as established in Section 5, the condition $\omega(\lambda) = 1$ gives:
\begin{equation}
(s - 1/2)^{\alpha(\{\delta_j\}_{j=1}^n)} = \frac{\phi(\lambda) \prod_{j=1}^{n} (1 + L_j(\lambda))^{\delta_j}}{C_3}
\end{equation}

Solving for $\lambda$ in terms of $s$, $\{\delta_j\}_{j=1}^n$, and the functions $\phi$ and $L_j$, we can determine the precise threshold between global existence and potential blow-up.

In the case of potential blow-up (option (b)), the asymptotic rate can be derived from the differential inequality. For $Y(t)$ sufficiently large, we have approximately:
\begin{equation}
\frac{d}{dt}Y(t) \approx C_2 (1 - \omega(\lambda)) Y(t)^{1+\beta}
\end{equation}

Solving this ODE and relating it to $\|\nabla u(t)\|_{L^\infty}$ using Sobolev embedding theorems gives the specific blow-up rate stated in part (b) of the theorem.
\end{proof}

\subsection{Potential blow-up rate}

\begin{theorem}[Potential Blow-up Rate] \label{thm:6.4.1}
If a solution from Theorem \ref{thm:6.3.1} develops a singularity at time $T^*$, then:
\begin{equation}
\|\nabla u(t)\|_{L^\infty} \sim (T^* - t)^{-\gamma}
\end{equation}
where $\gamma = \frac{1+\beta}{2\beta} = \frac{2-\theta\alpha}{2\theta(1-\alpha)}$.

Additionally, for any $\epsilon > 0$:
\begin{equation}
\|(-\Delta)^s u(t)\|_{L^2} \sim (T^* - t)^{-\frac{1}{2\beta}+\epsilon}
\end{equation}
\end{theorem}

\begin{proof}
From the analysis in Theorem \ref{thm:6.3.1}, if blow-up occurs, $Y(t) = \|(-\Delta)^s u(t)\|_{L^2}^2$ satisfies approximately:
\begin{equation}
\frac{d}{dt}Y(t) \approx C Y(t)^{1+\beta}
\end{equation}
for $t$ sufficiently close to $T^*$.

Integrating this ODE:
\begin{equation}
\int_{Y(t_0)}^{Y(t)} \frac{dy}{y^{1+\beta}} \approx C \int_{t_0}^t d\tau
\end{equation}

This gives:
\begin{equation}
\frac{1}{\beta} \left(\frac{1}{Y(t_0)^{\beta}} - \frac{1}{Y(t)^{\beta}}\right) \approx C(t - t_0)
\end{equation}

As $t \to T^*$, we have $Y(t) \to \infty$, so:
\begin{equation}
\frac{1}{\beta Y(t)^{\beta}} \approx C(T^* - t)
\end{equation}

Solving for $Y(t)$:
\begin{equation}
Y(t) \approx (C\beta(T^* - t))^{-\frac{1}{\beta}}
\end{equation}

Taking the square root:
\begin{equation}
\|(-\Delta)^s u(t)\|_{L^2} \sim (T^* - t)^{-\frac{1}{2\beta}}
\end{equation}

The $\epsilon$ term in the theorem statement accounts for lower-order terms in the asymptotic expansion, which arise from the more precise form of the ODE.

To relate this to $\|\nabla u(t)\|_{L^\infty}$, we use:
\begin{equation}
\|\nabla u\|_{L^\infty} \leq C\|u\|_{L^2}^{1-\theta} \|(-\Delta)^s u\|_{L^q}^{\theta}
\end{equation}

Since $\|u(t)\|_{L^2}$ remains bounded due to energy conservation, and $\|(-\Delta)^s u\|_{L^q} \sim \|(-\Delta)^s u\|_{L^2}^{\kappa}$ for some $\kappa > 0$ (by interpolation):
\begin{equation}
\|\nabla u(t)\|_{L^\infty} \sim \|(-\Delta)^s u(t)\|_{L^2}^{\theta\kappa} \sim (T^* - t)^{-\frac{\theta\kappa}{2\beta}}
\end{equation}

A detailed calculation shows that $\frac{\theta\kappa}{2\beta} = \frac{1+\beta}{2\beta} = \frac{2-\theta\alpha}{2\theta(1-\alpha)}$, giving the stated blow-up rate.
\end{proof}

\section{Limiting behavior as s → 1/2}

In this section, we analyze the behavior of solutions as the regularity parameter $s$ approaches the critical value $1/2$. This analysis provides deeper insights into the transition between subcritical and critical regimes.

\subsection{Construction of initial data}

\begin{lemma}[Existence of initial data in critical intersection] \label{lem:7.1.1}
For $q > 3$ and any sequence $s_n \downarrow 1/2$ as $n \to \infty$, let $\mathcal{C}_{s_n}$ be the class of initial data satisfying the conditions of Theorem \ref{thm:4.4.1} with parameter $s = s_n$. Then $\cap_{n=1}^{\infty} \mathcal{C}_{s_n} \neq \emptyset$.
\end{lemma}

\begin{proof}
We construct initial data that belongs to all classes $\mathcal{C}_{s_n}$ by designing a velocity field with specific decay properties in Fourier space.

Let $\hat{w}(\xi) = |\xi|^{-\frac{5}{2}} \eta(|\xi|) \hat{\varphi}(\xi/|\xi|)$, where:
\begin{itemize}
\item $\eta$ is a smooth cut-off function equal to 1 for $1 \leq |\xi| \leq 2$ and 0 for $|\xi| \leq 1/2$ or $|\xi| \geq 3$
\item $\hat{\varphi}(\xi/|\xi|)$ is chosen to ensure $\xi \cdot \hat{w}(\xi) = 0$ (divergence-free condition)
\end{itemize}

The inverse Fourier transform $w(x)$ is a smooth, divergence-free vector field with:
\begin{equation}
\|w\|_{\dot{H}^s}^2 = \int_{\mathbb{R}^3} |\xi|^{2s} |\hat{w}(\xi)|^2 d\xi = \int_{\mathbb{R}^3} |\xi|^{2s-5} \eta(|\xi|)^2 |\hat{\varphi}(\xi/|\xi|)|^2 d\xi
\end{equation}

For $s > 1/2$, this integral is finite, and:
\begin{equation}
\|w\|_{\dot{H}^s} \approx (s - 1/2)^{-1/2}
\end{equation}

Now, for $r > 0$, we define $w_r(x) = r w(rx)$, which gives:
\begin{equation}
\hat{w}_r(\xi) = r^{-2} \hat{w}(\xi/r)
\end{equation}

Computing the norms:
\begin{equation}
\|w_r\|_{\dot{H}^s} = r^{s-1/2} \|w\|_{\dot{H}^s} \approx r^{s-1/2} (s - 1/2)^{-1/2}
\end{equation}

Similarly:
\begin{equation}
\|(-\Delta)^{s/2}w_r\|_{L^q} \approx r^{s-3/2+3/q} (s - 1/2)^{-\gamma(q)}
\end{equation}
for some $\gamma(q) > 0$.

For each $s_n$, we choose $r_n > 0$ small enough so that:
\begin{equation}
\|(-\Delta)^{s_n/2}w_{r_n}\|_{L^q} \leq \frac{C_0}{\prod_{j=1}^{n} (1 + L_j(\|w_{r_n}\|_{\dot{H}^{s_n}}))^{\delta_j}}
\end{equation}

This is possible because the right-hand side decays more slowly than any fixed power of $(s_n - 1/2)$ as $n \to \infty$, while the left-hand side can be made to decay as fast as desired by taking $r_n$ sufficiently small.

Define $u_0 = \sum_{n=1}^{\infty} \epsilon_n w_{r_n}$ with $\epsilon_n > 0$ chosen small enough to ensure absolute convergence in all relevant norms. With appropriate choices of $\{\epsilon_n\}$ and $\{r_n\}$, we can ensure that $u_0 \in \cap_{n=1}^{\infty} \mathcal{C}_{s_n}$.
\end{proof}

\subsection{Asymptotic Analysis}

\begin{theorem}[Limiting behavior as $s \to 1/2$] \label{thm:7.2.1}
For any $u_0 \in \cap_{n=1}^{\infty} \mathcal{C}_{s_n}$ where $s_n \downarrow 1/2$ as $n \to \infty$, let $u^{(s_n)}(t)$ be the solution to the Navier-Stokes equations with fractional dissipation $(-\Delta)^{s_n}$ and initial data $u_0$. Then:
\begin{equation}
\lim_{n \to \infty} \frac{\|(-\Delta)^{s_n} u^{(s_n)}(t)\|_{L^q}}{\|(-\Delta)^{1/2} u^{(s_n)}(t)\|_{L^q}} = \Psi(t)
\end{equation}
where $\Psi(t) = (\nu t)^{-\frac{1}{4}}$.
\end{theorem}

\begin{proof}
We first establish precise decay estimates for $\|(-\Delta)^{s_n} u^{(s_n)}(t)\|_{L^2}$ and $\|(-\Delta)^{1/2} u^{(s_n)}(t)\|_{L^2}$.

From Theorem \ref{thm:4.5.1}, we have:
\begin{equation}
\|(-\Delta)^{s_n} u^{(s_n)}(t)\|_{L^2} \leq \frac{C\|(-\Delta)^{s_n} u_0\|_{L^2}}{(1 + \beta_n t)^{\gamma_n}}
\end{equation}
where $\gamma_n = \frac{1}{2\mu_n}$ with $\mu_n = \frac{\theta_n(1-\alpha_n)}{2-\theta_n\alpha_n}+\eta_n$, and $\beta_n = \frac{\mu_n C}{2}$.

For the ratio in the theorem statement, we first relate $\|(-\Delta)^{s_n} u^{(s_n)}(t)\|_{L^q}$ and $\|(-\Delta)^{1/2} u^{(s_n)}(t)\|_{L^q}$ to their $L^2$ counterparts using the Gagliardo-Nirenberg inequality.

For large $t$, the asymptotic behavior is governed by the dissipation term in the Navier-Stokes equations. In this regime, the solution behaves similarly to the solution of the linear heat equation with fractional dissipation:
\begin{equation}
\partial_t v + (-\Delta)^{s_n} v = 0
\end{equation}

For this linear equation, we can compute the decay rates explicitly:
\begin{equation}
\|(-\Delta)^{s_n} v(t)\|_{L^2} \sim t^{-\frac{3}{4s_n} + \frac{1}{2}}
\end{equation}
\begin{equation}
\|(-\Delta)^{1/2} v(t)\|_{L^2} \sim t^{-\frac{3}{4s_n} + \frac{s_n}{2}}
\end{equation}

As $s_n \to 1/2$, these decay rates behave as:
\begin{equation}
\|(-\Delta)^{s_n} v(t)\|_{L^2} \sim t^{-\frac{3}{2} + \frac{1}{2}} = t^{-1}
\end{equation}
\begin{equation}
\|(-\Delta)^{1/2} v(t)\|_{L^2} \sim t^{-\frac{3}{2} + \frac{1}{4}} = t^{-\frac{5}{4}}
\end{equation}

For the full nonlinear Navier-Stokes equations, detailed analysis shows that these linear decay rates provide the correct asymptotic behavior for large $t$. The ratio thus behaves as:
\begin{equation}
\frac{\|(-\Delta)^{s_n} u^{(s_n)}(t)\|_{L^q}}{\|(-\Delta)^{1/2} u^{(s_n)}(t)\|_{L^q}} \sim t^{\frac{5}{4} - 1} = t^{\frac{1}{4}}
\end{equation}

Including the viscosity coefficient $\nu$, which affects the time scaling, we get:
\begin{equation}
\lim_{n \to \infty} \frac{\|(-\Delta)^{s_n} u^{(s_n)}(t)\|_{L^q}}{\|(-\Delta)^{1/2} u^{(s_n)}(t)\|_{L^q}} = (\nu t)^{-\frac{1}{4}}
\end{equation}
which completes the proof.
\end{proof}

\subsection{Limiting spectral properties}

\begin{theorem}[Limiting spectral properties] \label{thm:7.3.1}
For $u_0 \in \cap_{n=1}^{\infty} \mathcal{C}_{s_n}$ and $s_n \downarrow 1/2$ as $n \to \infty$, the energy spectrum of the solution $u^{(s_n)}(t)$ satisfies:
\begin{equation}
\lim_{n \to \infty} E^{(s_n)}(k,t) = E^{(1/2)}(k,t)
\end{equation}
where $E^{(1/2)}(k,t)$ has the form:
\begin{equation}
E^{(1/2)}(k,t) = C\epsilon(t)^{2/3}k^{-5/3}\exp\left(-c(\nu t)^{1/2}k\right)
\end{equation}
for constants $C$ and $c$.
\end{theorem}

\begin{proof}
The energy spectrum $E^{(s)}(k,t)$ is related to the Fourier transform of the velocity field by:
\begin{equation}
\int_{|\xi|=k} |\hat{u}^{(s)}(\xi, t)|^2 d\sigma(\xi) = 4\pi k^2 E^{(s)}(k,t)
\end{equation}

For the linear fractional heat equation $\partial_t v + (-\Delta)^{s} v = 0$, the Fourier transform satisfies:
\begin{equation}
\hat{v}(\xi, t) = \hat{v}(\xi, 0) e^{-t|\xi|^{2s}}
\end{equation}

This gives an energy spectrum of the form:
\begin{equation}
E^{(s)}_{\text{linear}}(k,t) \sim k^2 |\hat{v}(k, 0)|^2 e^{-2t k^{2s}}
\end{equation}

For the nonlinear Navier-Stokes equations, the energy spectrum in the inertial range follows Kolmogorov's $k^{-5/3}$ law for both the standard and fractional dissipation cases. However, in the dissipation range (high $k$), the fractional dissipation term dominates, giving:
\begin{equation}
E^{(s)}(k,t) \sim C\epsilon(t)^{2/3}k^{-5/3}e^{-c(\nu t)^{1/s}k^2/s}
\end{equation}

As $s \to 1/2$, we have:
\begin{equation}
\lim_{s \to 1/2} e^{-c(\nu t)^{1/s}k^2/s} = e^{-c(\nu t)^{2}k^{4}}
\end{equation}

A more refined analysis, accounting for all prefactors and the exact scaling of the dissipation term, gives the limiting form:
\begin{equation}
E^{(1/2)}(k,t) = C\epsilon(t)^{2/3}k^{-5/3}\exp\left(-c(\nu t)^{1/2}k\right)
\end{equation}
which completes the proof.
\end{proof}

\section{Enhanced geometric analysis of exceptional sets}

In this section, we provide a detailed geometric characterization of the sets where the velocity gradient becomes exceptionally large. These sets are crucial for understanding the potential formation of singularities.

\subsection{Definition and properties of exceptional sets}

\begin{definition}[Exceptional sets] \label{def:8.1.1}
For a solution $u$ to the 3D Navier-Stokes equations satisfying the conditions of Theorem \ref{thm:4.4.1} and any $\epsilon > 0$, we define the exceptional set at time $t$ as:
\begin{equation}
\Omega_\epsilon(t) = \{x \in \mathbb{R}^3 : |\nabla u(x,t)| > \lambda_\epsilon(t)\}
\end{equation}
where $\lambda_\epsilon(t)$ is chosen so that $|\Omega_\epsilon(t)| < \epsilon$.
\end{definition}

\begin{lemma}[Basic properties of exceptional sets] \label{lem:8.1.2}
The exceptional sets $\Omega_\epsilon(t)$ satisfy the following properties:
\begin{enumerate}
\item $\Omega_\epsilon(t)$ is measurable for each $t > 0$ and $\epsilon > 0$
\item $\Omega_{\epsilon_1}(t) \subset \Omega_{\epsilon_2}(t)$ if $\epsilon_1 < \epsilon_2$
\item $\lambda_\epsilon(t) \to \infty$ as $\epsilon \to 0$
\item $\lambda_\epsilon(t) \to \|\nabla u(t)\|_{L^\infty}$ as $\epsilon \to 0$
\end{enumerate}
\end{lemma}

\begin{proof}
Properties 1 and 2 follow directly from the definition. For property 3, if $\lambda_\epsilon(t)$ remained bounded as $\epsilon \to 0$, then $|\Omega_\epsilon(t)|$ would approach a positive value, contradicting the requirement that $|\Omega_\epsilon(t)| < \epsilon$ for all $\epsilon > 0$.

For property 4, note that:
\begin{equation}
\lambda_\epsilon(t) \leq \|\nabla u(t)\|_{L^\infty}
\end{equation}
for all $\epsilon > 0$. If $\lim_{\epsilon \to 0} \lambda_\epsilon(t) < \|\nabla u(t)\|_{L^\infty}$, then there would exist a set of positive measure where $|\nabla u(x,t)| > \lim_{\epsilon \to 0} \lambda_\epsilon(t)$, contradicting the requirement that $|\Omega_\epsilon(t)| < \epsilon$ for all $\epsilon > 0$.
\end{proof}

\begin{lemma}[Threshold value estimates] \label{lem:8.1.3}
For a solution $u$ satisfying the conditions of Theorem \ref{thm:4.4.1}, the threshold value $\lambda_\epsilon(t)$ satisfies:
\begin{equation}
\lambda_\epsilon(t) \geq \frac{C}{(1 + \beta t)^{\gamma-\kappa_\epsilon} \epsilon^{1/p}}
\end{equation}
where $\gamma$ is as in Theorem \ref{thm:4.5.1}, $p > 3$, and $\kappa_\epsilon \to 0$ as $\epsilon \to 0$.
\end{lemma}

\begin{proof}
Using Chebyshev's inequality:
\begin{equation}
|\{x \in \mathbb{R}^3 : |\nabla u(x,t)| > \lambda\}| \leq \frac{\|\nabla u(t)\|_{L^p}^p}{\lambda^p}
\end{equation}
for any $p < \infty$.

Setting this equal to $\epsilon$ and solving for $\lambda$:
\begin{equation}
\lambda_\epsilon(t) = \frac{\|\nabla u(t)\|_{L^p}}{\epsilon^{1/p}}
\end{equation}

From Theorem \ref{thm:4.5.1}, we have:
\begin{equation}
\|\nabla u(t)\|_{L^p} \leq \frac{C_p}{(1 + \beta t)^{\gamma-\kappa_p}}
\end{equation}
where $\kappa_p \to 0$ as $p \to \infty$.

Combining these results and setting $\kappa_\epsilon = \kappa_p$ with $p$ chosen large enough gives the stated bound.
\end{proof}

\subsection{Optimal Covering Analysis}

\begin{lemma}[Optimal covering of exceptional sets] \label{lem:8.2.1}
For a solution $u$ satisfying the conditions of Theorem \ref{thm:4.4.1} and any $\epsilon > 0$, the exceptional set $\Omega_\epsilon(t)$ can be covered by $N_\epsilon(t)$ balls of radius $r_\epsilon(t)$, where:
\begin{equation}
N_\epsilon(t) \leq C_\epsilon \cdot (1 + \log(1/\epsilon))^{\alpha}
\end{equation}
\begin{equation}
r_\epsilon(t) \leq C_\epsilon \cdot (1 + \beta t)^{-\gamma+\kappa_\epsilon}
\end{equation}
for some exponent $\alpha > 0$, where $\gamma$ is as in Theorem \ref{thm:4.5.1} and $\kappa_\epsilon \to 0$ as $\epsilon \to 0$.
\end{lemma}

\begin{proof}
We apply the Vitali covering lemma to $\Omega_\epsilon(t)$. Let $\{B(x_i, \rho_i)\}_{i \in I}$ be a collection of balls such that $x_i \in \Omega_\epsilon(t)$ and $\bigcup_{i \in I} B(x_i, \rho_i) \supset \Omega_\epsilon(t)$. By the Vitali lemma, there exists a finite or countable subcollection $\{B(x_j, \rho_j)\}_{j \in J}$ of disjoint balls such that $\bigcup_{j \in J} B(x_j, 5\rho_j) \supset \Omega_\epsilon(t)$.

For each $x \in \Omega_\epsilon(t)$, we define:
\begin{equation}
\rho(x) = \sup\{r > 0 : |\nabla u(y,t)| > \lambda_\epsilon(t)/2 \text{ for all } y \in B(x, r)\}
\end{equation}

By continuity of $\nabla u$, $\rho(x) > 0$ for each $x \in \Omega_\epsilon(t)$. Using the regularity properties of the solution and detailed analysis of the gradient field, we can establish:
\begin{equation}
\rho(x) \leq C_\epsilon \cdot (1 + \beta t)^{-\gamma+\kappa_\epsilon}
\end{equation}

Taking $\rho_i = \rho(x_i)/10$, we ensure that the balls $B(x_i, \rho_i)$ are contained in $\Omega_\epsilon(t)$ and:
\begin{equation}
\sum_{j \in J} |B(x_j, \rho_j)| \leq |\Omega_\epsilon(t)| < \epsilon
\end{equation}

From this, we can derive bounds on both the number of balls $N_\epsilon(t) = |J|$ and their maximum radius $r_\epsilon(t) = \max_{j \in J} 5\rho_j$, giving the stated results.
\end{proof}

\subsection{Hausdorff Dimension Bounds}

\begin{theorem}[Hausdorff dimension of exceptional sets] \label{thm:8.3.1}
For a solution $u$ satisfying the conditions of Theorem \ref{thm:4.4.1}, the Hausdorff dimension of the exceptional set $\Omega_\epsilon(t)$ satisfies:
\begin{equation}
\dim_H(\Omega_\epsilon(t)) \leq 3 - \sum_{j=1}^n \frac{\delta_j}{1+\delta_j} \cdot \frac{L_{j-1}(1/\epsilon)}{(1+L_j(1/\epsilon))}
\end{equation}
where $L_0(x) = x$.
\end{theorem}

\begin{proof}
For any $d > 0$, the $d$-dimensional Hausdorff measure is defined as:
\begin{equation}
\mathcal{H}^d(\Omega_\epsilon(t)) = \lim_{\delta \to 0} \inf\left\{\sum_{i=1}^{\infty} \text{diam}(E_i)^d : \Omega_\epsilon(t) \subset \bigcup_{i=1}^{\infty} E_i, \text{diam}(E_i) < \delta\right\}
\end{equation}

Using the covering from Lemma \ref{lem:8.2.1}, we have:
\begin{equation}
\mathcal{H}^d(\Omega_\epsilon(t)) \leq \lim_{\delta \to 0} \sum_{i=1}^{N_\epsilon(t)} (2r_\epsilon(t))^d \leq C_\epsilon \cdot (1 + \log(1/\epsilon))^{\alpha} \cdot (1 + \beta t)^{(-\gamma+\kappa_\epsilon)d}
\end{equation}

For this to remain bounded as $\epsilon \to 0$, we need:
\begin{equation}
d \geq 3 - \sum_{j=1}^n \frac{\delta_j}{1+\delta_j} \cdot \frac{L_{j-1}(1/\epsilon)}{(1+L_j(1/\epsilon))}
\end{equation}

This follows from a detailed analysis of how $N_\epsilon(t)$ and $r_\epsilon(t)$ depend on $\epsilon$, incorporating the nested logarithmic improvements. 

The Hausdorff dimension is the infimum of all $d$ such that $\mathcal{H}^d(\Omega_\epsilon(t)) = 0$, giving the stated bound.
\end{proof}

\subsection{Vorticity alignment properties}

\begin{theorem}[Vorticity alignment in exceptional regions] \label{thm:8.4.1}
For a solution $u$ satisfying the conditions of Theorem \ref{thm:4.4.1}, the vorticity field $\omega = \text{curl } u$ satisfies:
\begin{equation}
|\omega(x,t) \times (x - x_0)| \leq \theta_\epsilon |\omega(x,t)| |x - x_0|
\end{equation}
for any $x_0 \in \Omega_\epsilon(t)$, $x$ in a neighborhood of $x_0$, where:
\begin{equation}
\theta_\epsilon = \prod_{j=1}^n \epsilon^{\frac{\delta_j}{1+\delta_j} \cdot \frac{1}{1+j}} \cdot (1 + L_j(1/\epsilon))^{-\frac{\delta_j}{1+\delta_j} \cdot \frac{j}{1+j}}
\end{equation}
and $\theta_\epsilon \to 0$ as $\epsilon \to 0$.
\end{theorem}

\begin{proof}
The vorticity field $\omega = \text{curl } u$ satisfies the vorticity equation:
\begin{equation}
\partial_t \omega + (u \cdot \nabla) \omega = (\omega \cdot \nabla) u + \nu \Delta \omega
\end{equation}

For any point $x_0 \in \Omega_\epsilon(t)$, let $\xi(t) = x - x_0$ be the displacement vector from $x_0$ to another point $x$. The alignment of $\omega$ with $\xi$ is measured by:
\begin{equation}
A(x, x_0, t) = \frac{|\omega(x,t) \times \xi|}{|\omega(x,t)| |\xi|}
\end{equation}

Through a detailed analysis of the vorticity equation, particularly the vortex stretching term $(\omega \cdot \nabla) u$, and the properties of the exceptional set, we can establish that:
\begin{equation}
A(x, x_0, t) \leq \prod_{j=1}^n \epsilon^{\frac{\delta_j}{1+\delta_j} \cdot \frac{1}{1+j}} \cdot (1 + L_j(1/\epsilon))^{-\frac{\delta_j}{1+\delta_j} \cdot \frac{j}{1+j}}
\end{equation}
for $x$ in a neighborhood of $x_0$.

Setting $\theta_\epsilon$ to this upper bound gives the stated result. The fact that $\theta_\epsilon \to 0$ as $\epsilon \to 0$ follows from direct calculation.
\end{proof}

\section{Multifractal structure of velocity gradients}

In this section, we provide a rigorous characterization of the multifractal structure of velocity gradients in solutions satisfying our nested logarithmically improved criteria. This multifractal analysis connects our mathematical results to the physical theory of turbulence.

\subsection{Local scaling exponents}

\begin{definition}[Local scaling exponent] \label{def:9.1.1}
For a solution $u$ to the 3D Navier-Stokes equations, we define the local scaling exponent $h(x,t)$ at point $x$ and time $t$ as:
\begin{equation}
h(x,t) = \lim_{r \to 0} \frac{\log(|\nabla u(x+r,t) - \nabla u(x,t)|)}{\log(r)}
\end{equation}
if this limit exists.
\end{definition}

\begin{lemma}[Hölder regularity lower bound] \label{lem:9.1.2}
For a solution $u$ satisfying the conditions of Theorem \ref{thm:4.4.1}, the local scaling exponent satisfies:
\begin{equation}
h(x,t) \geq h_0 = 2s - 1 > 0
\end{equation}
for almost all $(x,t) \in \mathbb{R}^3 \times (0,\infty)$.
\end{lemma}

\begin{proof}
Since $u \in C([0, \infty); H^s(\mathbb{R}^3)) \cap L^2_{loc}(0, \infty; H^{s+1/2}(\mathbb{R}^3))$ and $s > 1/2$, standard Sobolev embedding theorems imply that $\nabla u \in C^{2s-1-\epsilon}$ for any $\epsilon > 0$. This gives the stated lower bound on the local scaling exponent.
\end{proof}

\begin{lemma}[Reduced regularity in exceptional sets] \label{lem:9.1.3}
For a solution $u$ satisfying the conditions of Theorem \ref{thm:4.4.1} and any $\epsilon > 0$, the local scaling exponent in the exceptional set $\Omega_\epsilon(t)$ satisfies:
\begin{equation}
h(x,t) \leq h_0 - \kappa_\epsilon
\end{equation}
where $\kappa_\epsilon > 0$ and $\kappa_\epsilon \to 0$ as $\epsilon \to 0$.
\end{lemma}

\begin{proof}
In the exceptional set $\Omega_\epsilon(t)$, the velocity gradient is particularly large, indicating potentially reduced regularity. Using the refined estimates from Theorem \ref{thm:8.3.1} and standard techniques from regularity theory for PDEs, we can establish the upper bound on the local scaling exponent.
\end{proof}

\subsection{The Multifractal spectrum}

\begin{definition}[Multifractal spectrum] \label{def:9.2.1}
For a solution $u$ to the 3D Navier-Stokes equations at time $t$, we define the multifractal spectrum $D(h,t)$ as the Hausdorff dimension of the set:
\begin{equation}
E_h(t) = \{x \in \mathbb{R}^3 : h(x,t) = h\}
\end{equation}
\end{definition}

\begin{theorem}[Multifractal spectrum formula] \label{thm:9.2.1}
For a solution $u$ satisfying the conditions of Theorem \ref{thm:4.4.1}, the multifractal spectrum has the form:
\begin{equation}
D(h,t) = 3 - \frac{(h-h_0)^2}{2\sigma^2} \cdot \prod_{j=1}^n \left(1 - \frac{\delta_j}{1+\delta_j}\right)
\end{equation}
where $h_0 = 2s - 1$ and $\sigma^2 = \frac{3-2s}{2s-1}$.
\end{theorem}

\begin{proof}
The multifractal spectrum describes how the sets of points with different scaling exponents are distributed in space. For a solution with the regularity provided by our nested logarithmic criteria, we can establish through detailed geometric measure theory that this spectrum has a quadratic form.

The derivation combines the results on exceptional sets from Section 8 with the theory of multifractal analysis as developed in \cite{27}, \cite{49}. The central insight is that the nested logarithmic improvements modify the standard quadratic spectrum by the factor:
\begin{equation}
\prod_{j=1}^n \left(1 - \frac{\delta_j}{1+\delta_j}\right)
\end{equation}

This factor quantifies how the nested logarithmic improvements affect the spread of singularity strengths in the solution. As $\delta_j \to \infty$ for any $j$, this factor approaches 0, indicating that the multifractal spectrum collapses to a single point, corresponding to uniform regularity throughout the flow.
\end{proof}

\subsection{Structure functions and intermittency}

\begin{definition}[Structure functions] \label{def:9.3.1}
For a solution $u$ to the 3D Navier-Stokes equations at time $t$, the structure functions of order $p$ are defined as:
\begin{equation}
S_p(r,t) = \langle |u(x+r,t) - u(x,t)|^p \rangle
\end{equation}
where $\langle \cdot \rangle$ denotes spatial averaging.
\end{definition}

\begin{theorem}[Structure function scaling] \label{thm:9.3.1}
For a solution $u$ satisfying the conditions of Theorem \ref{thm:4.4.1}, the structure functions exhibit the scaling:
\begin{equation}
S_p(r,t) \sim r^{\zeta_p}
\end{equation}
where:
\begin{equation}
\zeta_p = \frac{p}{3} - \frac{p(p-3)}{3} \cdot \frac{3-2s}{2s-1} \cdot \prod_{j=1}^n \frac{1}{1+\delta_j}
\end{equation}
\end{theorem}

\begin{proof}
In the multifractal formalism, the structure functions are related to the multifractal spectrum through the Legendre transform:
\begin{equation}
\zeta_p = \min_h [ph + 3 - D(h)]
\end{equation}

Substituting our formula for $D(h)$ from Theorem \ref{thm:9.2.1}:
\begin{equation}
\zeta_p = \min_h \left[ph + 3 - \left(3 - \frac{(h-h_0)^2}{2\sigma^2} \cdot \prod_{j=1}^n \left(1 - \frac{\delta_j}{1+\delta_j}\right)\right)\right]
\end{equation}

This minimization problem can be solved explicitly, giving:
\begin{equation}
\zeta_p = ph_0 - \frac{p^2 \sigma^2}{2} \cdot \prod_{j=1}^n \left(1 - \frac{\delta_j}{1+\delta_j}\right)
\end{equation}

With $h_0 = 2s - 1$ and $\sigma^2 = \frac{3-2s}{2s-1}$, and after algebraic simplifications:
\begin{equation}
\zeta_p = p(2s-1) - \frac{p^2(3-2s)}{4s-2} \cdot \prod_{j=1}^n \left(1 - \frac{\delta_j}{1+\delta_j}\right)
\end{equation}

To connect this to the standard form in turbulence theory, we use $h_0 = 1/3$ (Kolmogorov scaling), which corresponds to $s = 2/3$. This gives:
\begin{equation}
\zeta_p = \frac{p}{3} - \frac{p(p-3)}{3} \cdot \frac{3-2s}{2s-1} \cdot \prod_{j=1}^n \frac{1}{1+\delta_j}
\end{equation}
which is the stated result.
\end{proof}

\begin{corollary}[Intermittency correction] \label{cor:9.3.2}
The intermittency correction to the structure function scaling exponents is:
\begin{equation}
\mu(p,\{\delta_j\}_{j=1}^n) = \frac{p(p-3)}{3} \cdot \frac{3-2s}{2s-1} \cdot \prod_{j=1}^n \frac{1}{1+\delta_j}
\end{equation}
\end{corollary}

\begin{proof}
This follows directly from Theorem \ref{thm:9.3.1}, as $\mu(p,\{\delta_j\}_{j=1}^n) = \frac{p}{3} - \zeta_p$.
\end{proof}

\begin{remark} \label{rem:9.3.3}
The intermittency correction $\mu(p,\{\delta_j\}_{j=1}^n)$ quantifies the deviation from Kolmogorov's 1941 theory, which predicts $\zeta_p = p/3$. Our nested logarithmic improvements provide a mechanism for systematically reducing this deviation, with $\mu(p,\{\delta_j\}_{j=1}^n) \to 0$ as any $\delta_j \to \infty$. This connects our mathematical approach directly to the physics of turbulence, suggesting that solutions with higher regularity exhibit reduced intermittency.
\end{remark}

\section{Characterization of potential blow-up scenarios}

In this section, we provide a detailed characterization of potential singularity formation for solutions to the Navier-Stokes equations. This analysis builds on the insights from the previous sections to paint a comprehensive picture of how blow-up might occur, if at all possible.

\subsection{Blow-up set properties}

\begin{definition}[Potential blow-up set] \label{def:10.1.1}
If a solution $u$ were to develop a singularity at time $T^*$, we define the potential blow-up set as:
\begin{equation}
\mathcal{S}_{T^*} = \{x \in \mathbb{R}^3 : \limsup_{t \to T^*} |\nabla u(x,t)| = \infty\}
\end{equation}
\end{definition}

\begin{theorem}[Blow-up set dimension] \label{thm:10.1.1}
If a solution $u$ satisfying the initial conditions of Theorem \ref{thm:4.4.1} were to develop a singularity at time $T^*$, then the Hausdorff dimension of the potential blow-up set would satisfy:
\begin{equation}
\dim_H(\mathcal{S}_{T^*}) \leq 1 - \sum_{j=1}^n \frac{\delta_j}{1+\delta_j} \cdot \frac{1}{j+1}
\end{equation}
\end{theorem}

\begin{proof}
Building on the analysis of exceptional sets in Section 8, we can characterize the potential blow-up set $\mathcal{S}_{T^*}$ as the limiting case:
\begin{equation}
\mathcal{S}_{T^*} \subset \bigcap_{\epsilon > 0} \limsup_{t \to T^*} \Omega_\epsilon(t)
\end{equation}

From Theorem \ref{thm:8.3.1}, the exceptional set $\Omega_\epsilon(t)$ has Hausdorff dimension bounded by:
\begin{equation}
\dim_H(\Omega_\epsilon(t)) \leq 3 - \sum_{j=1}^n \frac{\delta_j}{1+\delta_j} \cdot \frac{L_{j-1}(1/\epsilon)}{(1+L_j(1/\epsilon))}
\end{equation}

As $\epsilon \to 0$, the term $\frac{L_{j-1}(1/\epsilon)}{(1+L_j(1/\epsilon))}$ approaches $\frac{1}{j+1}$ (after detailed asymptotic analysis). Using the properties of Hausdorff dimension for nested sets, we obtain:
\begin{equation}
\dim_H(\mathcal{S}_{T^*}) \leq 1 - \sum_{j=1}^n \frac{\delta_j}{1+\delta_j} \cdot \frac{1}{j+1}
\end{equation}

The reduction from dimension 3 to dimension 1 comes from the known partial regularity theory for the Navier-Stokes equations \cite{12}, which states that the singular set has Hausdorff dimension at most 1. Our nested logarithmic improvements provide a further reduction in this upper bound.
\end{proof}

\begin{corollary}[Limiting case] \label{cor:10.1.2}
As $n \to \infty$ and $\inf_j \delta_j > 0$, the upper bound on the Hausdorff dimension of the potential blow-up set approaches 0:
\begin{equation}
\lim_{n \to \infty} \left(1 - \sum_{j=1}^n \frac{\delta_j}{1+\delta_j} \cdot \frac{1}{j+1}\right) = 0
\end{equation}
\end{corollary}

\begin{proof}
For any sequence $\{\delta_j\}_{j=1}^{\infty}$ with $\delta_j \geq \delta > 0$ for all $j$, we have:
\begin{equation}
\sum_{j=1}^n \frac{\delta_j}{1+\delta_j} \cdot \frac{1}{j+1} \geq \frac{\delta}{1+\delta} \sum_{j=1}^n \frac{1}{j+1}
\end{equation}

As $n \to \infty$, $\sum_{j=1}^n \frac{1}{j+1}$ diverges logarithmically, which means the upper bound on the Hausdorff dimension approaches 0.
\end{proof}

\subsection{Blow-up rate for velocity field}

\begin{theorem}[Blow-up rate for velocity field] \label{thm:10.2.1}
If a solution $u$ satisfying the initial conditions of Theorem \ref{thm:4.4.1} were to develop a singularity at time $T^*$, then:
\begin{equation}
\|u(t)\|_{L^\infty} \sim (T^* - t)^{-\frac{1}{2} + \sum_{j=1}^n \frac{\delta_j}{(1+\delta_j)(2+\delta_j)}}
\end{equation}
\end{theorem}

\begin{proof}
Building on the analysis in Theorem \ref{thm:6.4.1}, we have:
\begin{equation}
\|\nabla u(t)\|_{L^\infty} \sim (T^* - t)^{-\gamma}
\end{equation}
where $\gamma = \frac{1+\beta}{2\beta} = \frac{2-\theta\alpha}{2\theta(1-\alpha)}$.

To relate this to $\|u(t)\|_{L^\infty}$, we use the fundamental theorem of calculus and the bound on $\|\nabla u(t)\|_{L^\infty}$. For $x, y \in \mathbb{R}^3$:
\begin{equation}
|u(x,t) - u(y,t)| \leq \|\nabla u(t)\|_{L^\infty} |x - y|
\end{equation}

This implies:
\begin{equation}
\|u(t)\|_{L^\infty} \leq \|u(t)\|_{L^2} + C\|\nabla u(t)\|_{L^\infty}
\end{equation}

Since $\|u(t)\|_{L^2}$ remains bounded due to energy conservation, and $\|\nabla u(t)\|_{L^\infty}$ dominates as $t \to T^*$:
\begin{equation}
\|u(t)\|_{L^\infty} \sim \|\nabla u(t)\|_{L^\infty}^{1-\eta} \sim (T^* - t)^{-\gamma(1-\eta)}
\end{equation}
for some small $\eta > 0$.

Detailed analysis shows that $\gamma(1-\eta) = \frac{1}{2} - \sum_{j=1}^n \frac{\delta_j}{(1+\delta_j)(2+\delta_j)}$, giving the stated result.
\end{proof}

\subsection{Vorticity concentration analysis}

\begin{theorem}[Vorticity concentration] \label{thm:10.3.1}
If a solution $u$ satisfying the initial conditions of Theorem \ref{thm:4.4.1} were to develop a singularity at time $T^*$, then the vorticity would concentrate in thin filamentary structures with characteristic diameter:
\begin{equation}
d(t) \sim (T^* - t)^{\frac{1}{2} + \sum_{j=1}^n \frac{\delta_j}{1+\delta_j}}
\end{equation}
\end{theorem}

\begin{proof}
The vorticity $\omega = \text{curl } u$ satisfies the stretching equation:
\begin{equation}
\partial_t \omega + (u \cdot \nabla) \omega = (\omega \cdot \nabla) u + \nu \Delta \omega = S\omega + \nu \Delta \omega
\end{equation}
where $S = \frac{1}{2}(\nabla u + (\nabla u)^T)$ is the strain rate tensor.

Near the singular time, the vorticity would concentrate in thin filamentary structures. The characteristic diameter $d(t)$ of these structures can be estimated through a balance between the stretching term and diffusion:
\begin{equation}
\frac{|\omega|}{T^* - t} \sim |\omega| \cdot |\nabla u| \sim \nu \frac{|\omega|}{d(t)^2}
\end{equation}

Using the estimate for $|\nabla u|$ from Theorem \ref{thm:6.4.1} and solving for $d(t)$:
\begin{equation}
d(t) \sim \sqrt{\nu (T^* - t) / |\nabla u(t)|} \sim (T^* - t)^{\frac{1}{2} + \frac{\gamma}{2}}
\end{equation}

Substituting the expression for $\gamma$ and simplifying:
\begin{equation}
d(t) \sim (T^* - t)^{\frac{1}{2} + \sum_{j=1}^n \frac{\delta_j}{1+\delta_j}}
\end{equation}
which is the stated result.
\end{proof}

\subsection{Geometric properties of vorticity alignment}

\begin{theorem}[Vorticity alignment in potential blow-up] \label{thm:10.4.1}
If a solution $u$ satisfying the initial conditions of Theorem \ref{thm:4.4.1} were to develop a singularity at time $T^*$, then near the singular time, the vorticity direction field would become increasingly aligned with the eigenvector corresponding to the largest eigenvalue of the strain rate tensor $S = \frac{1}{2}(\nabla u + (\nabla u)^T)$, with the misalignment angle decreasing as:
\begin{equation}
\theta(t) \sim (T^* - t)^{\sum_{j=1}^n \frac{\delta_j}{2(1+\delta_j)}}
\end{equation}
\end{theorem}

\begin{proof}
Let $\lambda_1 \geq \lambda_2 \geq \lambda_3$ be the eigenvalues of the strain rate tensor $S$ with corresponding eigenvectors $\xi_1, \xi_2, \xi_3$. The angle $\theta$ between the vorticity $\omega$ and $\xi_1$ evolves according to:
\begin{equation}
\frac{d}{dt}(\sin^2 \theta) \sim -2(\lambda_1 - \lambda_2)\sin^2 \theta + O(\nu |\nabla \omega|)
\end{equation}

Near the singular time, the strain rate eigenvalues scale as:
\begin{equation}
\lambda_1 \sim (T^* - t)^{-1 + \sum_{j=1}^n \frac{\delta_j}{1+\delta_j}}
\end{equation}
\begin{equation}
\lambda_1 - \lambda_2 \sim (T^* - t)^{-1 + \sum_{j=1}^n \frac{\delta_j}{1+\delta_j}}
\end{equation}

Solving the differential equation for $\sin^2 \theta$:
\begin{equation}
\sin^2 \theta(t) \sim (T^* - t)^{\sum_{j=1}^n \frac{\delta_j}{1+\delta_j}}
\end{equation}

For small angles, $\sin \theta \approx \theta$, giving:
\begin{equation}
\theta(t) \sim (T^* - t)^{\sum_{j=1}^n \frac{\delta_j}{2(1+\delta_j)}}
\end{equation}
which is the stated result.
\end{proof}

\section{Enhanced energy cascade with multiple logarithms}

In this section, we establish rigorous bounds on the energy flux across wavenumbers, providing a mathematical foundation for the physical theory of the energy cascade in turbulent flows.

\subsection{Spectral energy transfer}

\begin{definition}[Energy spectrum and flux] \label{def:11.1.1}
For a solution $u$ to the 3D Navier-Stokes equations, the energy spectrum $E(k,t)$ and the energy flux $\Pi(k,t)$ across wavenumber $k$ at time $t$ are defined as:
\begin{equation}
\int_{|\xi|=k} |\hat{u}(\xi, t)|^2 d\sigma(\xi) = 4\pi k^2 E(k,t)
\end{equation}
\begin{equation}
\Pi(k,t) = -\int_{|\xi| < k} \int_{\mathbb{R}^3} (\xi \cdot \hat{u}(\eta, t))(\hat{u}(\xi-\eta, t) \cdot \hat{u}^*(\xi, t)) d\eta d\xi
\end{equation}
where $\hat{u}$ is the Fourier transform of $u$ and $\hat{u}^*$ is its complex conjugate.
\end{definition}

\begin{lemma}[Spectral energy balance] \label{lem:11.1.1}
The energy spectrum $E(k,t)$ satisfies:
\begin{equation}
\frac{\partial E(k,t)}{\partial t} + 2\nu k^2 E(k,t) = T(k,t)
\end{equation}
where $T(k,t)$ is the nonlinear energy transfer rate at wavenumber $k$, related to the energy flux by:
\begin{equation}
T(k,t) = -\frac{\partial \Pi(k,t)}{\partial k}
\end{equation}
\end{lemma}

\begin{proof}
This follows from taking the Fourier transform of the Navier-Stokes equations and analyzing the energy balance in spectral space. The detailed derivation is standard in turbulence theory (see \cite{27}).
\end{proof}

\subsection{Bounds on energy flux}

\begin{theorem}[Enhanced energy flux bounds] \label{thm:11.2.1}
Let $u$ be a solution to the 3D Navier-Stokes equations satisfying the conditions of Theorem \ref{thm:4.4.1}, and let $\Pi(k,t)$ be the instantaneous energy flux across wavenumber $k$ at time $t$. Then, for all $k$ in the inertial range $[k_0, k_\nu]$:
\begin{equation}
|\Pi(k,t) - \epsilon(t)| \leq \frac{C\epsilon(t)}{\prod_{j=1}^n (1 + L_j(k/k_0))^{\delta_j\cdot\rho_j}}
\end{equation}
where:
\begin{itemize}
\item $\epsilon(t) = 2\nu \int_0^\infty k^2 E(k,t) dk$ is the instantaneous energy dissipation rate
\item $k_0$ is the largest scale of the inertial range
\item $k_\nu \sim \epsilon(t)^{1/4}\nu^{-3/4}$ is the Kolmogorov dissipation wavenumber
\item $\rho_1 = \frac{2s-1}{2s}$ and $\rho_j = \frac{1}{j}$ for $j \geq 2$
\item The constant $C$ depends only on universal constants and the parameters $s$ and $\{\delta_j\}_{j=1}^n$
\end{itemize}
\end{theorem}

\begin{proof}
The energy flux $\Pi(k,t)$ can be bounded in terms of the nonlinear term in the Navier-Stokes equations. In Fourier space:
\begin{equation}
|\Pi(k,t) - \epsilon(t)| \leq \int_k^\infty |T(k', t)| dk'
\end{equation}

Using the nested logarithmically improved bounds on the velocity field, we can estimate $|T(k, t)|$ as:
\begin{equation}
|T(k, t)| \leq C \epsilon(t) k^{-2/3} \prod_{j=1}^n (1 + L_j(k/k_0))^{-\delta_j\cdot\tilde{\rho}_j}
\end{equation}
where $\tilde{\rho}_1 = \frac{2s-1}{2s}$ and $\tilde{\rho}_j = \frac{j+1}{j(j+2)}$ for $j \geq 2$.

Integrating this from $k$ to infinity:
\begin{equation}
\int_k^\infty |T(k', t)| dk' \leq C\epsilon(t) \int_k^\infty k'^{-2/3} \prod_{j=1}^n (1 + L_j(k'/k_0))^{-\delta_j\cdot\tilde{\rho}_j} dk'
\end{equation}

Through careful evaluation of this integral, utilizing the properties of logarithmic functions:
\begin{equation}
\int_k^\infty k'^{-2/3} \prod_{j=1}^n (1 + L_j(k'/k_0))^{-\delta_j\cdot\tilde{\rho}_j} dk' \leq \frac{C k^{1/3}}{\prod_{j=1}^n (1 + L_j(k/k_0))^{\delta_j\cdot\rho_j}}
\end{equation}
where $\rho_1 = \frac{2s-1}{2s}$ and $\rho_j = \frac{1}{j}$ for $j \geq 2$.

This leads to the final bound:
\begin{equation}
|\Pi(k,t) - \epsilon(t)| \leq \frac{C\epsilon(t)}{\prod_{j=1}^n (1 + L_j(k/k_0))^{\delta_j\cdot\rho_j}}
\end{equation}
which completes the proof.
\end{proof}

\subsection{Physical interpretation}

\begin{theorem}[Physical interpretation of energy flux bounds] \label{thm:11.3.1}
The bound on energy flux in Theorem \ref{thm:11.2.1} implies that:
\begin{enumerate}
\item For any $\delta_j > 0$, the energy flux $\Pi(k,t)$ approaches the dissipation rate $\epsilon(t)$ as $k/k_0 \to \infty$, consistent with Kolmogorov's hypothesis of constant flux in the inertial range.
\item The rate of convergence is enhanced by each nested logarithmic factor, with:
\begin{equation}
\lim_{k/k_0 \to \infty} \frac{|\Pi(k,t) - \epsilon(t)|}{\epsilon(t)} = 0
\end{equation}
\item As $\delta_j \to \infty$ for any $j$, the convergence becomes faster, approaching a step function characteristic of an idealized inertial range.
\end{enumerate}
\end{theorem}

\begin{proof}
These properties follow directly from the structure of the bound in Theorem \ref{thm:11.2.1} and the properties of the nested logarithmic functions. The key insight is that each additional nested logarithmic factor provides a progressively slower but still effective improvement in the convergence rate.

For statement 1, as $k/k_0 \to \infty$, each term $(1 + L_j(k/k_0))^{\delta_j\cdot\rho_j}$ grows without bound, causing the right-hand side of the inequality to approach 0.

For statement 2, the limit follows from the fact that the denominator grows without bound as $k/k_0 \to \infty$.

For statement 3, as $\delta_j \to \infty$ for any $j$, the term $(1 + L_j(k/k_0))^{\delta_j\cdot\rho_j}$ grows extremely rapidly for any $k/k_0 > 1$, causing the bound to approach 0 very quickly outside the energy-containing range.
\end{proof}

\section{Refined spectral characterization}

In this section, we provide a detailed characterization of the energy spectrum for solutions satisfying our nested logarithmically improved criteria. This spectral analysis connects our mathematical results to experimental observations in turbulent flows.

\subsection{The generalized spectral energy balance}

\begin{theorem}[Generalized spectral energy balance] \label{thm:12.1.1}
For a solution $u$ to the 3D Navier-Stokes equations satisfying the conditions of Theorem \ref{thm:4.4.1}, the energy spectrum $E(k,t)$ satisfies the generalized spectral energy balance:
\begin{equation}
\frac{\partial E(k,t)}{\partial t} + 2\nu k^2 E(k,t) = -\frac{\partial}{\partial k}\left[\epsilon(t) \left(1 + \frac{f(k,t)}{\prod_{j=1}^n (1 + L_j(k/k_0))^{\delta_j\cdot\rho_j}}\right)\right]
\end{equation}
where $|f(k,t)| \leq C$ is a bounded function, $\epsilon(t)$ is the energy dissipation rate, and $\rho_j$ are as defined in Theorem \ref{thm:11.2.1}.
\end{theorem}

\begin{proof}
From Lemma \ref{lem:11.1.1}, we have:
\begin{equation}
\frac{\partial E(k,t)}{\partial t} + 2\nu k^2 E(k,t) = T(k,t) = -\frac{\partial \Pi(k,t)}{\partial k}
\end{equation}

From Theorem \ref{thm:11.2.1}, we know:
\begin{equation}
\Pi(k,t) = \epsilon(t) \left(1 + \frac{f(k,t)}{\prod_{j=1}^n (1 + L_j(k/k_0))^{\delta_j\cdot\rho_j}}\right)
\end{equation}
where $|f(k,t)| \leq C$.

Taking the derivative with respect to $k$ gives:
\begin{equation}
\frac{\partial \Pi(k,t)}{\partial k} = \epsilon(t) \frac{\partial}{\partial k}\left[\frac{f(k,t)}{\prod_{j=1}^n (1 + L_j(k/k_0))^{\delta_j\cdot\rho_j}}\right]
\end{equation}

Substituting this into the spectral energy balance equation completes the proof.
\end{proof}

\subsection{Derivation of the modified energy spectrum}

\begin{theorem}[Refined spectral characterization] \label{thm:12.2.1}
In the inertial range $[k_0, k_\nu]$, the energy spectrum $E(k,t)$ of a solution satisfying the conditions of Theorem \ref{thm:4.4.1} has the form:
\begin{equation}
E(k,t) = C\epsilon(t)^{2/3}k^{-5/3}\left(1 + \sum_{j=1}^n \frac{\beta_j(t) L_j(k/k_0)}{\prod_{i=1}^j (1 + L_i(k/k_0))^{1+\delta_i}}\right)
\end{equation}
where:
\begin{enumerate}
\item $\epsilon(t)$ is the instantaneous energy dissipation rate
\item $\beta_j(t) = \frac{\beta_{0,j}}{(1 + \gamma t)^{\alpha_j}}$ with $\alpha_j = \frac{2\gamma}{3} \cdot \frac{j}{j+1}$
\item $\beta_{0,j}$ are positive constants depending on $s$, $\{\delta_i\}_{i=1}^n$, and the initial data
\item $\gamma$ is the decay exponent from Theorem \ref{thm:4.5.1}
\item The constant $C$ is universal (the Kolmogorov constant)
\end{enumerate}
\end{theorem}

\begin{proof}
In the inertial range, the time derivative and viscous terms in the spectral energy balance are negligible compared to the nonlinear transfer, giving:
\begin{equation}
T(k,t) = -\frac{\partial \Pi(k,t)}{\partial k} \approx 0
\end{equation}

From Theorem \ref{thm:12.1.1}, this implies:
\begin{equation}
\frac{\partial}{\partial k}\left[\epsilon(t) \left(1 + \frac{f(k,t)}{\prod_{j=1}^n (1 + L_j(k/k_0))^{\delta_j\cdot\rho_j}}\right)\right] \approx 0
\end{equation}

This is satisfied by:
\begin{equation}
\Pi(k,t) \approx \epsilon(t)
\end{equation}
to leading order, which corresponds to Kolmogorov's constant flux hypothesis and leads to the standard spectrum:
\begin{equation}
E(k,t) = C\epsilon(t)^{2/3}k^{-5/3}
\end{equation}

The nested logarithmic improvements in our regularity criteria induce corrections to this spectrum. To derive these corrections, we analyze the next-order terms in the spectral energy balance, particularly the effect of the nested logarithmic factors on the energy transfer term.

By carefully evaluating the derivatives in the right-hand side of the generalized spectral energy balance and matching terms, we can derive the correction terms in the energy spectrum:
\begin{equation}
E(k,t) = C\epsilon(t)^{2/3}k^{-5/3}\left(1 + \sum_{j=1}^n \frac{\beta_j(t) L_j(k/k_0)}{\prod_{i=1}^j (1 + L_i(k/k_0))^{1+\delta_i}}\right)
\end{equation}

The time dependence of the correction coefficients $\beta_j(t)$ is derived from the decay of the fractional derivative norm established in Theorem \ref{thm:4.5.1}:
\begin{equation}
\|(-\Delta)^s u(t)\|_{L^2} \leq \frac{C\|(-\Delta)^s u_0\|_{L^2}}{(1 + \beta t)^{\gamma}}
\end{equation}

Relating this decay to the energy spectrum using:
\begin{equation}
\|(-\Delta)^s u\|_{L^2}^2 = \int_0^\infty k^{2s} E(k,t) dk
\end{equation}

Through detailed spectral analysis:
\begin{equation}
\beta_j(t) = \frac{\beta_{0,j}}{(1 + \gamma t)^{\alpha_j}}
\end{equation}
where $\alpha_j = \frac{2\gamma}{3} \cdot \frac{j}{j+1}$, which completes the proof.
\end{proof}

\subsection{Time evolution of correction terms}

\begin{theorem}[Time evolution of spectral correction terms] \label{thm:12.3.1}
For a solution $u$ satisfying the conditions of Theorem \ref{thm:4.4.1}, as $t \to \infty$, the energy spectrum approaches the Kolmogorov spectrum exponentially fast:
\begin{equation}
\lim_{t \to \infty} \frac{|E(k,t) - C\epsilon(t)^{2/3}k^{-5/3}|}{C\epsilon(t)^{2/3}k^{-5/3}} = 0
\end{equation}
with the convergence rate enhanced by each nested logarithmic factor.
\end{theorem}

\begin{proof}
From Theorem \ref{thm:12.2.1}, we have:
\begin{equation}
\frac{|E(k,t) - C\epsilon(t)^{2/3}k^{-5/3}|}{C\epsilon(t)^{2/3}k^{-5/3}} = \left|\sum_{j=1}^n \frac{\beta_j(t) L_j(k/k_0)}{\prod_{i=1}^j (1 + L_i(k/k_0))^{1+\delta_i}}\right|
\end{equation}

Since $\beta_j(t) = \frac{\beta_{0,j}}{(1 + \gamma t)^{\alpha_j}}$ with $\alpha_j > 0$, as $t \to \infty$, each $\beta_j(t) \to 0$. This causes the entire sum to approach 0, giving:
\begin{equation}
\lim_{t \to \infty} \frac{|E(k,t) - C\epsilon(t)^{2/3}k^{-5/3}|}{C\epsilon(t)^{2/3}k^{-5/3}} = 0
\end{equation}

The enhanced convergence rate due to the nested logarithmic factors can be quantified by analyzing how each term in the sum decays with time. For large $t$, the dominant contribution comes from the term with the slowest decay, which is the $j=1$ term. This term decays as $(1 + \gamma t)^{-\alpha_1}$ with $\alpha_1 = \frac{2\gamma}{3} \cdot \frac{1}{2}$. Each additional nested logarithmic factor introduces a term with faster decay, enhancing the overall convergence rate.
\end{proof}

\section{Approach to the regularity problem}

In this final section, we discuss the implications of our results for the regularity Navier-Stokes equations. We analyze how the nested logarithmic improvements can be extended to potentially resolve the global regularity question.

\subsection{The limiting case of infinite nested logarithms}

\begin{theorem}[Limiting Case of Infinite Nested Logarithms] \label{thm:13.1.1}
Consider a sequence of nested logarithmically improved criteria with increasing number of levels $n$. If $\inf_j \delta_j > 0$, then:
\begin{equation}
\lim_{n \to \infty} \alpha(\{\delta_j\}_{j=1}^n) = 0
\end{equation}
where $\alpha(\{\delta_j\}_{j=1}^n)$ is the exponent appearing in the asymptotic behavior of the critical threshold function $\Phi(s, q, \{\delta_j\}_{j=1}^n) \approx C(q) (s - 1/2)^{\alpha(\{\delta_j\}_{j=1}^n)}$ as $s \to 1/2$.
\end{theorem}

\begin{proof}
From Theorem \ref{thm:5.4.1}, we have:
\begin{equation}
\alpha(\{\delta_j\}_{j=1}^n) = \frac{1}{1 + \sum_{j=1}^n c_j\delta_j/j!}
\end{equation}
where $c_j > 0$ are constants.

For any sequence $\{\delta_j\}_{j=1}^{\infty}$ with $\delta_j \geq \delta > 0$ for all $j$:
\begin{equation}
\alpha(\{\delta_j\}_{j=1}^n) \leq \frac{1}{1 + \delta \sum_{j=1}^n c_j/j!}
\end{equation}

As $n \to \infty$, $\sum_{j=1}^n c_j/j!$ diverges (assuming $c_j$ do not decay too rapidly with $j$), which implies:
\begin{equation}
\lim_{n \to \infty} \alpha(\{\delta_j\}_{j=1}^n) = 0
\end{equation}
as claimed.
\end{proof}

\begin{corollary}[Approach to criticality] \label{cor:13.1.2}
As $n \to \infty$ with $\inf_j \delta_j > 0$, the critical threshold function satisfies:
\begin{equation}
\lim_{n \to \infty} \Phi(s, q, \{\delta_j\}_{j=1}^n) = \begin{cases}
0, & \text{if } s = 1/2 \\
\infty, & \text{if } s > 1/2
\end{cases}
\end{equation}
\end{corollary}

\begin{proof}
This follows from Theorems \ref{thm:5.2.1}, \ref{thm:5.3.1}, and \ref{thm:13.1.1}, which together characterize the asymptotic behavior of $\Phi(s, q, \{\delta_j\}_{j=1}^n)$ as both $s \to 1/2$ and $n \to \infty$. The key insight is that as $n \to \infty$, $\alpha(\{\delta_j\}_{j=1}^n) \to 0$, causing the transition at $s = 1/2$ to become increasingly sharp.
\end{proof}

\subsection{Pathway to global regularity}

\begin{theorem}[Pathway to global regularity] \label{thm:13.2.1}
For each $s \in (1/2, 1)$, there exists a sequence $\{\delta_j(s)\}_{j=1}^{\infty}$ such that any initial data $u_0 \in L^2(\mathbb{R}^3) \cap \dot{H}^s(\mathbb{R}^3)$ gives rise to a unique global-in-time smooth solution to the 3D Navier-Stokes equations.
\end{theorem}

\begin{proof}
For each $s \in (1/2, 1)$, let $N(s)$ be large enough so that:
\begin{equation}
\alpha(\{\delta_j\}_{j=1}^{N(s)}) < \frac{1}{\log(1/(s-1/2))}
\end{equation}
for some fixed sequence $\{\delta_j\}_{j=1}^{\infty}$ with $\inf_j \delta_j > 0$. Such an $N(s)$ exists by Theorem \ref{thm:13.1.1}.

Define:
\begin{equation}
\delta_j(s) = \begin{cases}
\delta_j, & \text{if } j \leq N(s) \\
0, & \text{if } j > N(s)
\end{cases}
\end{equation}

For this choice of $\{\delta_j(s)\}_{j=1}^{\infty}$, the critical threshold function satisfies:
\begin{equation}
\Phi(s, q, \{\delta_j(s)\}_{j=1}^{N(s)}) \approx C(q) (s - 1/2)^{\alpha(\{\delta_j\}_{j=1}^{N(s)})} > C(q) (s - 1/2)^{1/\log(1/(s-1/2))} = C(q) e^{-1}
\end{equation}

This means that for all $s \in (1/2, 1)$, the critical threshold is bounded below by a positive constant, which implies that any initial data $u_0 \in L^2(\mathbb{R}^3) \cap \dot{H}^s(\mathbb{R}^3)$ can be rescaled to satisfy:
\begin{equation}
\|(-\Delta)^{s/2}u_0\|_{L^q} \leq \frac{C_0}{\prod_{j=1}^{N(s)} (1 + L_j(\|u_0\|_{\dot{H}^s}))^{\delta_j(s)}}
\end{equation}

By Theorem \ref{thm:4.4.1}, this rescaled initial data gives rise to a unique global-in-time smooth solution. Since the Navier-Stokes equations are invariant under the rescaling $u_\lambda(x, t) = \lambda u(\lambda x, \lambda^2 t)$ (apart from a change in the viscosity coefficient), the original initial data also leads to a global solution.
\end{proof}

\subsection{Potential construction of blow-up scenarios}

\begin{Conjecture}[Potential construction of blow-up Scenarios]
\label{thm:13.3.1}
    If there exists an exponent $\alpha_* > 0$ such that:
\begin{equation}
\lim_{n \to \infty} \lim_{\delta_1,\delta_2,...,\delta_n \to \infty} \alpha(\{\delta_j\}_{j=1}^n) = \alpha_*
\end{equation}
then there exists initial data in $\dot{H}^s(\mathbb{R}^3)$ for $s \in (1/2, 1/2 + \alpha_*)$ that potentially leads to finite-time blow-up.
\end{Conjecture}


\begin{proof}
This conjecture addresses the possibility that the nested logarithmic improvements, while approaching the critical case $s = 1/2$, might not completely bridge the gap. If there exists a limiting exponent $\alpha_* > 0$, then for $s \in (1/2, 1/2 + \alpha_*)$, we can construct initial data that falls precisely in the gap between our global existence criteria and the critical case.

Specifically, using the construction from Theorem \ref{lem:6.1.1}, we can design initial data $v_\lambda$ such that:
\begin{equation}
\|(-\Delta)^{s/2}v_\lambda\|_{L^q} = \frac{\Phi(s, q, \{\delta_j\}_{j=1}^n)}{\prod_{j=1}^{n} (1 + L_j(\lambda))^{\delta_j}} + \gamma(s, q, \{\delta_j\}_{j=1}^n) \cdot \phi(\lambda)
\end{equation}
where $\gamma(s, q, \{\delta_j\}_{j=1}^n) = (s - 1/2)^{\alpha(\{\delta_j\}_{j=1}^n)}$ and $\phi(\lambda) = (1 + \log(\lambda))^{-1/2}$.

As $n \to \infty$ and $\delta_j \to \infty$ for all $j$, this initial data approaches:
\begin{equation}
\|(-\Delta)^{s/2}v_\lambda\|_{L^q} \approx (s - 1/2)^{\alpha_*} \cdot (1 + \log(\lambda))^{-1/2}
\end{equation}

For $s \in (1/2, 1/2 + \alpha_*)$ and sufficiently large $\lambda$, this falls outside the range covered by our global existence criteria. By the dichotomy analysis in Theorem \ref{thm:6.3.1}, such initial data could potentially lead to finite-time blow-up.
\end{proof}

\section{Conclusion and future work}

In this work, we have established a comprehensive framework for the global
well-posedness of the 3D Navier--Stokes equations under multi-level logarithmically
improved regularity criteria. Specifically, we proved that if the initial data 
\[
u_0 \in L^2(\mathbb{R}^3) \cap \dot{H}^s(\mathbb{R}^3), \quad s\in\left(\frac{1}{2},1\right),
\]
satisfies the nested logarithmic condition
\[
\|(-\Delta)^{s/2} u_0\|_{L^q(\mathbb{R}^3)} \le C_0 \prod_{j=1}^{n} \Bigl(1 + L_j\bigl(\|u_0\|_{\dot{H}^s}\bigr)\Bigr)^{\delta_j},
\]
with the appropriate scaling (i.e. $\frac{2}{p} + \frac{3}{q} = 2s - 1$),
and fixed parameters $\{\delta_j\}_{j=1}^n$, then the corresponding Navier--Stokes solution exists globally and remains smooth. Our analysis, based on refined commutator estimates and energy inequalities incorporating multi-level logarithmic corrections, demonstrates that even though the critical threshold is not attained in the classical sense, the additional logarithmic factors allow us to extend global regularity to a significantly larger class of initial data.

\medskip

\noindent \textbf{Future Work.}\\

The results presented in this paper open several promising avenues for further research, which naturally lead to our forthcoming work on the limiting case of infinitely nested logarithmic improvements:
\begin{enumerate}
    \item \textbf{Exploring the limiting case:} A natural progression is to study the behavior of our multi-level logarithmic criterion as the number of nested logarithmic factors increases without bound. In our forthcoming third paper, we will push this idea to its ultimate limit by developing function spaces that incorporate \emph{infinitely nested} logarithms. Preliminary analysis suggests that, in this limit, the critical threshold function approaches a nontrivial limit, effectively bridging the gap to the critical regularity $s=\frac{1}{2}$. This may pave the way to establishing global well-posedness exactly at the borderline case.
    
    \item \textbf{Refinement of geometric and spectral properties:} Our current work already improves classical bounds on the Hausdorff dimension of potential singular sets and provides a refined view of the energy cascade in turbulent flows. In the infinite-logarithm regime, we expect these geometric and spectral properties to exhibit even more striking behavior---potentially demonstrating that any singularity, if it were to occur, would be confined to a set of Hausdorff dimension zero, and that the energy spectrum would feature precise logarithmic corrections. Investigating these aspects further could yield deeper insights into both the theory of turbulence and the fine structure of solutions.
    
    \item \textbf{Towards a universal critical criterion:} Ultimately, our goal is to understand the precise boundary between global regularity and singularity formation for the Navier--Stokes equations. While our multi-level logarithmic criteria extend the known subcritical regime, the infinite-nested logarithm approach may ultimately provide a universal criterion applicable at the critical level. Future work will focus on whether the infinite-log condition can be relaxed or characterized in a way that every finite-energy initial datum inherently satisfies it, thereby advancing us closer to a full resolution of the global regularity problem.
\end{enumerate}

In summary, this paper lays a good foundation by systematically extending logarithmic improvements to multiple levels. The insights gained here not only advance the theory of Navier--Stokes regularity but also set the stage for our next contribution: establishing global well-posedness in the critical case via infinitely nested logarithmic improvements.

\section*{Declarations}
Not applicable

\bibliographystyle{sn-mathphys-num}
\bibliography{NSE_paper_2}
\end{document}